\documentclass[preprint,11pt,a4paper,nopreprintline]{elsarticle}
\usepackage[utf8]{inputenc}
\usepackage[T1]{fontenc}
\usepackage{lmodern}
\usepackage[english]{babel}

\usepackage{amsmath,amssymb,amsthm,mathtools,mathrsfs}
\usepackage{booktabs,array,graphicx,xcolor}
\usepackage[unicode,hidelinks]{hyperref}
\newcommand{\R}{\mathbb R}
\newcommand{\diver}{\operatorname{div}}
\newcommand{\Ric}{\operatorname{Ric}}

\newcommand{\tr}{\operatorname{tr}}
\newcommand{\vol}{\operatorname{vol}}

\newcommand{\dd}{\,\mathrm d}

\newtheorem{theorem}{Theorem}
\newtheorem{lemma}{Lemma}
\newtheorem{proposition}{Proposition}
\newtheorem{corollary}{Corollary}
\theoremstyle{definition}

\newtheorem{example}{Example}

\newcommand{\Rm}{\operatorname{Rm}}

\journal{Journal of Mathematical Analysis and Applications}
\biboptions{sort&compress}

\newcommand{\manuscriptfigure}[1]{%
  \IfFileExists{#1}{\includegraphics[width=.90\linewidth]{#1}}{}%
}

\begin{document}
\begin{frontmatter}

\title{Dissipative Geometric Flows and Certified Phase-Volume Contraction}
\author[ipme]{I.B. Furtat}
\ead{cainenash@mail.ru}
\affiliation[ipme]{organization={IPME RAS},city={Saint Petersburg},country={Russia}}

\begin{abstract}
This paper develops an approach to uniform phase-volume contraction for prescribed dynamical systems. The main result is that every positive rate strictly below the invariant-measure bound can be certified in finite search time by a regularized Sobolev gradient flow. The resulting positive density ensures exponential contraction of the weighted phase volume with unit prefactor. For each fixed regularization parameter, the flow converges exponentially to a unique minimizer. As the regularization vanishes, these minimizers converge to the feasible potential of minimum Sobolev norm for a fixed feasibility margin.
The certification procedure rests on a H\"older estimate that converts a bound on the integral penalty into a pointwise bound on the weighted divergence over the entire invariant domain. This estimate yields a stopping criterion for the search. Finite termination is also established in growing finite-dimensional trial spaces for a staged scheme with exact gradient evaluations and controlled positive quadrature. The scheme requires neither a supplied feasibility margin nor a bound on the norm of a feasible potential. A disk example shows that the invariant-measure bound can be approached by smooth certificates even when no $C^1$ certificate attains it.
Analytical properties of novel local D-Liouville, D-Yamabe, and D-Ricci evolutions driven by divergence penalties are established. For the unnormalized D-Yamabe flow, a scaling obstruction to exponential stability of stationary conformal factors is established. For the D-Ricci flow, the transformation law for the prescribed vector field under DeTurck pullback is derived. This law specifies how the vector field must transform when the gauge term is removed. Examples on the torus and sphere identify obstructions to smooth convergence of the geometric evolutions.
\end{abstract}

\end{frontmatter}

\section{Introduction}

The evolution of a set of initial states provides information about a dynamical system that is difficult to obtain from individual trajectories. Its volume measures how much phase space the set occupies as it moves. A uniform contraction estimate quantifies how quickly this volume decreases, including during transients. Such an estimate can hold in a suitably weighted volume even when the divergence in the reference volume changes sign. This leads to two questions: which contraction rates are feasible, and how can a density that certifies a prescribed rate with unit prefactor be constructed and verified throughout the invariant domain?

Density criteria have been used to establish almost-global convergence \cite{Rantzer01,PrajnaParriloRantzer2004}, while differential contraction provides bounds on the separation of nearby trajectories \cite{Lohmiller98,Forni14}. For weighted phase-volume contraction, earlier work provides a geometric interpretation and a search over a fixed log-affine family \cite{Furtat26_Flow,Furtat26_Density}. Auxiliary-function duality gives the sharp rate bound over smooth densities, and finite-time averaging yields feasible potentials for rates strictly below that bound \cite{TobascoGoluskinDoering2018}. The questions addressed here are whether a search can be guaranteed to terminate with a certificate when only the target rate is supplied, how to select a feasible potential of minimum Sobolev norm, and whether the limiting rate itself is attained.

The link between weighted volume and its underlying geometry also motivates the study of evolution equations for the density or metric. Classical Yamabe and Ricci flows are driven by curvature \cite{Hamilton1982,Ye1994,Brendle07}. Weighted Yamabe flow uses the scalar curvature of a metric measure space \cite{Yan22}, and volume-preserving Willmore flow imposes an enclosed-volume constraint on a moving hypersurface \cite{Chen2026Willmore}. In the present setting, the geometry evolves in response to the weighted divergence of a prescribed vector field. This coupling raises two distinct issues. Constant rescaling leaves the divergence unchanged but alters the curvature contribution to the conformal evolution. A DeTurck pullback transforms both the metric and the prescribed field, which determines how the divergence term is carried into the ungauged equation.

The paper makes the following contributions to the certification of phase-volume contraction and the analysis of divergence-driven geometric flows:

(1) A finite-termination guarantee is proved for a Sobolev search in growing trial spaces. This extends the constructive scope of the sharp invariant-measure bound and of searches over fixed log-affine families \cite{TobascoGoluskinDoering2018,Furtat26_Density}. A H\"older estimate converts the integral penalty into a pointwise stopping test. With exact gradient evaluations and controlled positive quadrature, the search certifies any positive rate strictly below the bound using the target rate alone. At a fixed feasible margin, vanishing regularization selects the feasible potential of minimum Sobolev norm.

(2) A scalar log-density certifies volume decay for every material set in the invariant domain. The quantity controlled in \cite{GieslHafsteinMehrabinezhad2021} is transverse distance near an exponentially stable periodic orbit. The annular examples provide independently verified certificates for one-dimensional and polynomial--Fourier searches. A disk example establishes nonattainment at the limiting rate in $C^1$, even though smooth certificates approach that rate.

(3) The D-Liouville, D-Yamabe, and D-Ricci flows evolve, respectively, a positive density, a conformal factor, and a Riemannian metric, with each evolution incorporating a penalty on the weighted divergence of a prescribed vector field. Convergence of local D-Liouville evolution is established under regularization or directional coercivity. For unnormalized D-Yamabe flow, the different scaling laws of curvature and weighted divergence yield a precise alternative: a stationary factor lies in a family of constant rescalings or is linearly unstable \cite{Ye1994,Brendle07,Yan22}. Stable mean-zero states are constructed, and convergence for analytic penalties under uniform regularity bounds is proved using a \L{}ojasiewicz--Simon inequality \cite{WuGrasselliZheng2007,FeehanMaridakis2017}. For D-Ricci flow, the penalty tensor and its pullback law are derived. Removing the DeTurck term also transforms the prescribed field \cite{DeTurck1983}. With the field held fixed, a torus example gives different orders for the first nonconformal term in the gauged and ungauged evolutions, conditional on smooth solvability of the latter.

Sections~\ref{sec:prelim}--\ref{s:search2} establish feasibility and certification. 
Sections~\ref{sec:convergence}--\ref{sec:transport_limits} analyze the geometric flows and their implications for physical transport.

\section{Geometric Setting and Search Objective}\label{sec:prelim}

Let $(M,g)$ be a compact oriented Riemannian $n$-manifold, possibly with smooth boundary. Metrics and vector fields are smooth unless stated otherwise. Inward or tangent $f$ generates a forward semiflow $\varphi_t$, outward transport is considered only until exit. Write $dV_g$ for volume, $dS_g$ for boundary measure, $\nu$ for the outward unit normal, $f(v)=dv[f]=\langle f,\nabla_gv\rangle_g$, and $\Delta_g=\diver_g\nabla_g$ (nonpositive spectrum). On tensors, $\Delta_g=\nabla^i\nabla_i$ for the Levi-Civita connection $\nabla$. In coordinates, $\diver_gf=\partial_if^i+\Gamma^i_{ik}f^k$ and $\tr_gh=g^{ij}h_{ij}$, repeated indices are summed. Metric inequalities are inequalities of quadratic forms. For a volume form $m$, divergence is defined by $\mathcal L_fm=(\diver_mf)m$, where $\mathcal L$ denotes Lie derivative, $\Phi^*$ denotes pullback. The scalar curvature is $R_g=\tr_g\Ric(g)$, and $\Rm(g)$ denotes the Riemann tensor.

Write $C^{k,\beta}$, $H^s$, and $h^{k+\beta}$ for H\"older, Sobolev, and little H\"older spaces, with $h^{k+\beta}$ the closure of smooth functions or sections in $C^{k,\beta}$. Unless stated otherwise, $L^2$ uses $dV_g$. We write $a_+=\max(a,0)$, $\|\cdot\|_\infty$ for the uniform norm, $U\Subset W$ for compact containment, $\sigma(L)$ for the spectrum, and $DF$ for the Fr\'echet derivative. For a covector $\xi$, $f\cdot\xi=\xi(f)$, and $O(\delta^p)$ is measured in the stated norm.

Consider the system
\begin{equation}\label{eq:sys}
 \dot x=f(x),\qquad x\in M.
\end{equation}
All weights used in the geometric searches are strictly positive. For $w=e^\psi$ and $\alpha\in\mathbb R$, define the weighted divergence~\cite{Furtat26_Flow} by
\begin{equation}\label{eq:weighted_div}
 D_{w^\alpha}f=\diver_g f+\alpha f(\ln w).
\end{equation}
For $\alpha=1$, write $D_wf$. The associated conformal metric is
$\widetilde g_\alpha=w^{2\alpha/n}g$, $dV_{\widetilde g_\alpha}=w^\alpha dV_g$.
For $\alpha=0$, set $\widetilde g_0=g$ and $D_{w^0}f=\diver_gf$. Divergence depends only on the volume form. In particular, if $dV_h=\rho\,dV_g$ with smooth $\rho>0$, then
\begin{equation}\label{eq:volume_identity}
 \diver_hf=\diver_{\rho^{2/n}g}f
 =\diver_gf+f(\ln\rho),
 \qquad \diver_{\widetilde g_\alpha}f=D_{w^\alpha}f.
\end{equation}
The volume-form identity implies that smooth positive densities and Riemannian metrics realize the same divergences. Metric evolution allows additional geometric constraints. For Jordan-measurable $D$, write $\vol_h(D)=\int_DdV_h$.

On a compact oriented manifold without boundary, Stokes' theorem gives
\begin{equation}\label{eq:stokes_obstruction}
 \int_M\diver_hf\,dV_h=0.
\end{equation}
Consequently, uniformly negative divergence cannot hold on the entire manifold in any smooth metric. On such a manifold, the target must be imposed on a proper positively invariant domain.

On a closed manifold, integration of $\diver_g(uX)$ gives
$\int_M X(u)\,dV_g=-\int_M u\,\diver_g X\,dV_g$.

The static feasibility condition determines which contraction rates may be used as targets for an evolution. Let $K\subset M$ be a nonempty compact positively invariant set on which the forward flow exists for all time, and let $U$ be an open neighborhood of $K$. Write $b=\diver_g f$ and let $\mathscr I_f(K)$ denote the invariant probability measures on $K$. Define
\begin{equation}\label{eq:static_rate}
 \beta_K(f)=\max_{\mu\in\mathscr I_f(K)}\int_K b\,d\mu,
 \qquad \lambda_*(K,f)=-\beta_K(f).
\end{equation}
The following form of auxiliary-function duality~\cite{TobascoGoluskinDoering2018} provides the static bound and a smooth feasible competitor.
\begin{proposition}\label{prop:static_feasibility}
Under the preceding assumptions,
\begin{equation}\label{eq:static_duality}
 \beta_K(f)=\inf_{u\in C^\infty(U)}\max_{x\in K}\bigl(b(x)+f(u)(x)\bigr).
\end{equation}
Consequently, for every $0<\gamma<\lambda_*(K,f)$ there is a smooth potential $u$ on a neighborhood of $K$ such that $\diver_{e^u dV_g}f\le-\gamma$ on $K$. For $\gamma>\lambda_*(K,f)$ no $C^1$ potential satisfies this inequality. The value in~\eqref{eq:static_rate} is independent of the reference metric. Thus a positive contraction rate is feasible for some smooth density if and only if $\beta_K(f)<0$.
\end{proposition}
\begin{proof}
Invariant probability measures exist by compactness. Invariance gives $\int_Kf(u)\,d\mu=0$, so $\beta_K(f)\le\max_K(b+f(u))$ for every smooth $u$. For the reverse bound, set
$q_T(x)=\frac1T\int_0^T b(\varphi_t x)\,dt$, $u_T(x)=\frac1T\int_0^T(T-t)b(\varphi_t x)\,dt$.
For every invariant measure, $\int_Kq_T\,d\mu=\int_Kb\,d\mu$, hence $\max_Kq_T\ge\beta_K(f)$. Conversely, if $T_j\to\infty$ and $x_j$ maximizes $q_{T_j}$, any weak limit of the empirical measures $T_j^{-1}\int_0^{T_j}\delta_{\varphi_t x_j}\,dt$ is invariant. Here $\delta_x$ denotes unit Dirac mass. Indeed, shifting the integration interval by a fixed $s$ changes the integral of a continuous test function $v$ by at most $2s\|v\|_\infty/T_j$. A subsequence realizing the upper limit therefore gives $\limsup_{T\to\infty}\max_Kq_T\le\beta_K(f)$. Thus $\max_Kq_T\to\beta_K(f)$.

Differentiation along the flow and integration by parts give $b+f(u_T)=q_T$. For each finite $T$, smooth dependence on initial data makes $u_T$ smooth on a neighborhood of $K$. A cutoff equal to one near $K$ extends it smoothly to $U$ without changing this identity on $K$. Taking $T\to\infty$ proves~\eqref{eq:static_duality}, and sufficiently large finite $T$ gives feasibility for $0<\gamma<\lambda_*$. Integrating any feasible inequality against invariant measures gives the necessary bound $\gamma\le-\beta_K(f)$. Finally, if $dV_h=e^q dV_g$, invariance implies that the integrals of $\diver_hf=b+f(q)$ and $b$ agree. This proves metric independence.
\end{proof}
A strict bound on $K$ extends to a neighborhood by continuity. An all-time transport estimate requires positive invariance of the domain where the bound holds. Curvature constraints or prescribed boundary data can restrict feasibility. Total volume can be prescribed by adding a constant to $u$, which leaves $b+f(u)$ unchanged.

For a Jordan-measurable domain $D$ with piecewise smooth boundary, write
$\Omega_{w^\alpha}(D)=\int_D w^\alpha\,dV_g$.
If the metric and weight depend smoothly on physical time, the Liouville--Reynolds transport identity~\cite{Liouville,Reynolds,Furtat26_Flow} reads
\begin{equation}\label{eq:general_transport_th}
 \frac{d}{dt}\Omega_{w^\alpha}(\varphi_t(D),t)
 =\int_{\varphi_t(D)}\left(
 D_{w^\alpha}f+\alpha\partial_t\ln w
 +\frac12\tr_{g(t)}(\partial_tg(t))\right)w^\alpha\,dV_{g(t)}.
\end{equation}
Compute $D_{w^\alpha}f$ using $g(t)$ and apply the identity while the transported domain remains in the region of definition. Exponential $w^\alpha$-dissipativity means $\Omega_{w^\alpha}(\varphi_t(D),t)\le C e^{-\lambda t}\Omega_{w^\alpha}(D,0)$ for $C,\lambda>0$ independent of $D,t$. At $\lambda=0$ this is $w^\alpha$-dissipativity. Bounding the entire integrand in parentheses in~\eqref{eq:general_transport_th} by $-\lambda$ gives $C=1$.
 For $\lambda=0$, the term $w^\alpha$-dissipativity is used. If the expression in parentheses in~\eqref{eq:general_transport_th} is at most $-\lambda$, this estimate holds with $C=1$.

Once a smooth positive density or a Riemannian metric $h$ has been frozen, a sufficient pointwise condition for volume contraction on a positively invariant domain $U$ is $\diver_h f\le-\lambda$. In that case,
\begin{equation}\label{eq:vol_exp_D}
 \vol_h(\varphi_t(D))\le e^{-\lambda t}\vol_h(D),
 \qquad D\subset U,\quad t\ge0.
\end{equation}
Gronwall's inequality gives this bound from the volume-form evolution. Taking $dV_h=w^\alpha dV_g$ gives the weighted-volume estimate for ensembles of trajectories.
On a compact invariant domain, \eqref{eq:vol_exp_D} transfers from the frozen density $e^u$ to the background measure with prefactor $e^{\max u-\min u}$. Conversely, localization of a background estimate for every material set with fixed prefactor $C$ gives $\int_0^t b(\varphi_s x)\,ds\le\log C-\lambda t$. Integrating against invariant measures and letting $t\to\infty$ gives $\lambda\le\lambda_*$. Thus the supremal rate with a free prefactor is density-independent. The search instead constructs a pointwise bound for decay with unit prefactor in physical time.

%%%%%%%%%%%%%%%%%%%%%%%%%%%%%%%%%%%%

\section{Certified Sobolev Density Search}
\label{s:search2}

The search is formulated as a Sobolev gradient flow and is generally nonlocal. The regularity afforded by the Sobolev norm allows a bound on the integral penalty to be converted into a pointwise bound on the weighted divergence.
Let $\Omega$ be a compact connected domain with smooth boundary in an oriented Riemannian manifold $(M,g)$. The field $f$ is smooth on a neighborhood of $\Omega$ and points inward or tangent to its boundary, so $\Omega$ is positively invariant. Write $m=dV_g$, $b=\diver_gf$, and $L_fu=f(u)=du[f]$. For a positive density $\rho=e^u$,
$D_u=\diver_{e^um}f=b+L_fu.$
If $D_u\le-\lambda$ on $\Omega$, then every transported domain $E\subset\Omega$ satisfies
\begin{equation}
 \int_{\varphi_t(E)}e^u\dd m\le e^{-\lambda t}\int_Ee^u\dd m.\label{eq:transportstatic2}
\end{equation}
Physical time is denoted by $t$, while density-search time is denoted below by $\tau$. These variables are distinct.

For the compact invariant domain $\Omega$, write $\beta_\Omega=\beta_\Omega(f)$ and $\lambda_*=-\beta_\Omega$, as in~\eqref{eq:static_rate}. Proposition~\ref{prop:static_feasibility} provides a smooth feasible potential for every $0<\gamma<\lambda_*$. A positive target rate is possible precisely when $\beta_\Omega<0$.
Choose a search space in which convergence of the potential yields uniform control of its derivative along the vector field. Fix an integer $s>n/2+1$ and the Hilbert space $H=H^s(\Omega)$ with its usual full norm. The embedding theorem gives a compact embedding of $H$ into $C^1(\Omega)$. The operator $L_f:H\to L^2(\Omega,m)$ is bounded. Denote its Hilbert adjoint by $L_f^*$, using the $H$ inner product on the domain and the $L^2$ inner product on the codomain. The adjoint therefore includes the Riesz map for the chosen Sobolev inner product.

For a prescribed $\gamma>0$, define
\begin{equation}
 J_\gamma(u)=\frac12\|(b+L_fu+\gamma)_+\|_{L^2}^2,
 ~~~ E_\varepsilon(u)=J_\gamma(u)+\frac\varepsilon2\|u\|_H^2,
 ~~~\varepsilon>0.\label{eq:energy2}
\end{equation}
The positive part is taken pointwise. The functional is $C^{1,1}(H)$, although the scalar penalty has no continuous second derivative. Consider the flow in the $H$ metric:
\begin{equation}
 \partial_\tau u=-G_\varepsilon(u),\qquad
 G_\varepsilon(u)=L_f^*(b+L_fu+\gamma)_++\varepsilon u.
 \label{eq:Hflow2}
\end{equation}
This choice of inner product produces a generally nonlocal evolution.

\begin{theorem}
\label{t:main2}
Let $\mathcal C_\gamma=\{u\in H:b+L_fu\le-\gamma\text{ on }\Omega\}$ be nonempty.

1) For every $\varepsilon>0$, there is a unique minimizer $u_\varepsilon$ of $E_\varepsilon$. From every $u_0\in H$ the solution of~\eqref{eq:Hflow2} exists globally and satisfies
\begin{align}
 \|u(\tau)-u_\varepsilon\|_H
 &\le e^{-\varepsilon\tau}\|u_0-u_\varepsilon\|_H,\label{eq:rateH2}\\
 E_\varepsilon(u(\tau))-E_\varepsilon(u_\varepsilon)
 &\le e^{-2\varepsilon\tau}
 [E_\varepsilon(u_0)-E_\varepsilon(u_\varepsilon)].\label{eq:rateE2}
\end{align}

2) The minimum-norm potential $u^\dagger\in\mathcal C_\gamma$ is unique. As $\varepsilon\downarrow0$, $u_\varepsilon\to u^\dagger$ strongly in $H$ and $C^1$, and $J_\gamma(u_\varepsilon)\le(\varepsilon/2)\|u^\dagger\|_H^2$.

3) For every $0<\lambda<\gamma$, all sufficiently small fixed $\varepsilon$ give $D_{u(\tau)}\le-\lambda$ for all sufficiently large finite $\tau$. Freezing $e^{u(\tau)}$ gives~\eqref{eq:transportstatic2}. A smooth positive density can be used at any strictly smaller rate. In particular, if $\lambda_*>0$, choosing $\lambda<\gamma<\lambda_*$ gives finite-time attainment for every $0<\lambda<\lambda_*$.

4) Suppose $m(B_r(x)\cap\Omega)\ge\kappa_0 r^n$ for all $x\in\Omega$ and $0<r\le r_0$. Given $\eta>0$ and a bound $B>0$ on the $C^{0,\theta}$ seminorm of $R_u=(D_u+\gamma)_+$, with $0<\theta<1$, set
\begin{equation}\label{eq:stop2}
 r_\eta=\min\{r_0,(\eta/(2B))^{1/\theta}\},\qquad Q=\frac{\kappa_0\eta^2 r_\eta^n}{8}.
\end{equation}
If $J_\gamma(u)<Q$, then $D_u<-\gamma+\eta$ everywhere. For $s>n/2+1+\theta$, a bound $B$ is available from the Sobolev norm. In particular, if $\|v\|_H\le R$ for some $v\in\mathcal C_\gamma$, one may compute $B$ using the uniform bound $\|u_0\|_H+2R$ on the orbit. With this $B$ and $\eta=\gamma-\lambda$, sufficient choices for the target rate $\lambda$ are
\begin{equation}\label{eq:searchtime2}
 \varepsilon R^2\le Q,\qquad
 \tau>\frac{1}{2\varepsilon}
 \log\!\left(\max\!\left\{1,\frac{2E_\varepsilon(u_0)}{Q}\right\}\right).
\end{equation}

5) Let $H_N\subset H$ be nested finite-dimensional spaces of smooth functions with dense union. For fixed $\varepsilon$, the restricted minimizers $u_{\varepsilon,N}$ converge in $H$ to $u_\varepsilon$, and the restricted gradient flow satisfies the same rates as in~\eqref{eq:rateH2}--\eqref{eq:rateE2}. Given $\lambda<\gamma$, for sufficiently small fixed $\varepsilon$ and all sufficiently large $N$, the restricted flow reaches $D_u\le-\lambda$ in finite time from every $u_0\in H_N$.

\end{theorem}
\begin{proof}
The map $r\mapsto r_+$ is nonexpansive in $L^2$. Hence $G_\varepsilon$ is globally Lipschitz on $H$ with constant at most $\|L_f\|^2+\varepsilon$. Convexity of the penalty and the quadratic regularizer give
\begin{equation}
 \langle G_\varepsilon(u)-G_\varepsilon(v),u-v\rangle_H
 \ge\varepsilon\|u-v\|_H^2.\label{eq:mono2}
\end{equation}
The direct method in Hilbert space yields a unique minimum. Ordinary differential equation theory for globally Lipschitz vector fields on $H$ gives a global solution. Applying~\eqref{eq:mono2} to $u$ and $u_\varepsilon$ gives~\eqref{eq:rateH2}.

Along the flow, $E_\varepsilon'=-\|G_\varepsilon(u)\|_H^2$. Strong convexity gives
$E_\varepsilon(u)-E_\varepsilon(u_\varepsilon) \le\frac1{2\varepsilon}\|G_\varepsilon(u)\|_H^2$,
proving~\eqref{eq:rateE2}.

The set $\mathcal C_\gamma$ is convex and closed in $H$, so it has a unique metric projection of zero, $u^\dagger$. Comparison of the minima gives
\begin{equation}
 J_\gamma(u_\varepsilon)+\frac\varepsilon2\|u_\varepsilon\|_H^2
 \le\frac\varepsilon2\|u^\dagger\|_H^2.\label{eq:comparison2}
\end{equation}
In particular, $\|u_\varepsilon\|_H\le\|u^\dagger\|_H$. Every sequence $\varepsilon_j\downarrow0$ has a subsequence converging weakly in $H$ and strongly in $C^1$ to $v$. By~\eqref{eq:comparison2}, its penalty tends to zero. Thus $v\in\mathcal C_\gamma$, while weak lower semicontinuity of the norm gives $\|v\|_H\le\|u^\dagger\|_H$. Consequently, $v=u^\dagger$. Moreover,
$\|u^\dagger\|_H\le\liminf\|u_{\varepsilon_j}\|_H \le\limsup\|u_{\varepsilon_j}\|_H\le\|u^\dagger\|_H$.
Convergence of norms together with weak convergence gives strong convergence in $H$. Uniqueness of the limit removes the need to pass to a subsequence.

 Set $\eta=\gamma-\lambda>0$. For sufficiently small $\varepsilon$, $\|L_f(u_\varepsilon-u^\dagger)\|_\infty<\eta/2$. For all sufficiently large $\tau$, \eqref{eq:rateH2} gives $\|L_f(u(\tau)-u_\varepsilon)\|_\infty<\eta/2$. Their sum gives $D_{u(\tau)}<-\lambda$ on $\Omega$. Smooth approximation in $C^1$ preserves any smaller rate. Proposition~\ref{prop:static_feasibility} gives a smooth feasible potential for $\lambda<\gamma<\lambda_*$, with a strict margin if desired.

If $R_u(x_0)\ge\eta$, then $R_u\ge\eta/2$ on $B_{r_\eta}(x_0)\cap\Omega$. This implies $J_\gamma(u)\ge Q$. The positive-part map is $1$-Lipschitz, so the embedding $H^s\hookrightarrow C^{1,\theta}$ bounds $[R_u]_{C^{0,\theta}}$ in terms of $\|u\|_H$, $f$, $g$, and $\Omega$. A smooth compact domain has the stated lower volume bound. If the seminorm is zero, the residual is constant and can be checked directly. At a trial state the regularized energy also gives $\|u\|_H^2\le 2E_\varepsilon(u)/\varepsilon$. A spatial maximum computed with a guaranteed uniform error bound provides another test.

For a feasible $v$ with $\|v\|_H\le R$, comparison gives $\|u_\varepsilon\|_H\le R$ and $E_\varepsilon(u_\varepsilon)\le\varepsilon R^2/2$. The triangle inequality and~\eqref{eq:rateH2} give $\|u(\tau)\|_H\le\|u_0\|_H+2R$. By~\eqref{eq:rateE2},
$ J_\gamma(u(\tau))\le \frac{\varepsilon R^2}{2}
              + E_\varepsilon(u_0)e^{-2\varepsilon\tau}$.
Under~\eqref{eq:searchtime2} this is strictly below $Q$. One possible bound $R$ is the Sobolev norm of the finite-horizon potential in the proof of Proposition~\ref{prop:static_feasibility}, after choosing a horizon with the required divergence bound.

Choose $v_N\in H_N$ converging to $u_\varepsilon$ in $H$. Then
$ \tfrac{\varepsilon}{2}\|u_{\varepsilon,N}-u_\varepsilon\|_H^2
 \le E_\varepsilon(u_{\varepsilon,N})-E_\varepsilon(u_\varepsilon)
 \le E_\varepsilon(v_N)-E_\varepsilon(u_\varepsilon)\to0$.
The projected gradient obeys~\eqref{eq:mono2} on $H_N$, so the same distance and energy rates hold. For the sign guarantee, take $\varepsilon$ small enough that $\|L_f(u_\varepsilon-u^\dagger)\|_\infty<\eta/3$, then $N$ large enough that $\|L_f(u_{\varepsilon,N}-u_\varepsilon)\|_\infty<\eta/3$. The restricted distance rate makes the remaining error less than $\eta/3$ in finite time.
\end{proof}

Fix $0<\lambda<\gamma<\lambda_*$ and choose $0<\theta<1$ with $s>n/2+1+\theta$. Set $\varepsilon_j=2^{-j}$. At stage $j$ start from zero, run the full flow until $\tau_j=j/\varepsilon_j$, and test $J_\gamma<Q$ in~\eqref{eq:stop2} with $\eta=\gamma-\lambda$ and a H\"older bound computed from the current $H$ norm. At fixed $\varepsilon$, \eqref{eq:searchtime2} provides a time bound when a feasible norm bound is known. In the staged procedure a feasible competitor of norm at most $R$ exists by Proposition~\ref{prop:static_feasibility}, and $R$ is used only in the proof. The energy at the test is bounded by $\varepsilon_jR^2/2+e^{-2j}E_{\varepsilon_j}(0)$. Hence the penalty tends to zero and the $H$ norms remain bounded by $\|u\|_H^2\le2E_{\varepsilon_j}(u)/\varepsilon_j$. With a fixed Sobolev embedding bound the computed H\"older bounds remain bounded, so $Q$ stays bounded away from zero. The test therefore succeeds at a finite stage.

The same rule works for exact projected gradient steps on $H_j$, with fixed step $0<h\le(\|L_f\|^2+1)^{-1}$ and $k_j=\lceil j/(h\varepsilon_j)\rceil$ steps. A smooth competitor feasible for some $\gamma'\in(\gamma,\lambda_*)$ can be approximated in $H_j$ and hence in $C^1$. Thus the restricted problems have feasible competitors of uniformly bounded norm for all large $j$. The descent lemma and strong convexity give
$E_{\varepsilon_j}(u_k)-\min_{H_j}E_{\varepsilon_j}
\le (1-h\varepsilon_j)^k E_{\varepsilon_j}(0)$.
The preceding argument now uses $e^{-j}$ in place of $e^{-2j}$, with $e^{-j}/\varepsilon_j\to0$. For a computed potential, the stopping test requires a rigorous upper bound on its penalty, or an independent uniform upper bound on its divergence. In particular, an approximate quadrature value alone cannot certify the sign.

If only $\lambda<\lambda_*$ is supplied, set $\gamma_j=\lambda+1/j$ and use the same stage lengths. A single smooth competitor with margin $\gamma'\in(\lambda,\lambda_*)$, or its approximants on $H_j$, is feasible for all sufficiently large $j$. The initial penalties remain bounded and the stage energies are $O(2^{-j})$, so the Sobolev and H\"older bounds remain uniform. With $\eta_j=1/j$, the threshold~\eqref{eq:stop2} is bounded below by a positive constant times $j^{-(2+n/\theta)}$. The same test therefore succeeds at a finite stage for both the continuous and exact projected schemes.

Let $m_j$ be finite positive quadrature measures converging weakly to $m$. In the projected scheme replace $J_{\gamma_j}$ by $J_j(u)=\tfrac12\int R_j(u)^2\,dm_j$, where $R_j(u)=(D_u+\gamma_j)_+$. Use the same step count with $h$ chosen from a common Lipschitz bound. Such a bound is $1+\sup_j m_j(\Omega)\|L_f\|_{H\to C^0}^2$ for the regularized gradients. The preceding competitors are pointwise feasible, hence have zero quadrature penalty, so $J_j(u_j)=O(2^{-j})$ at the stage outputs $u_j$ and $\|u_j\|_H$ is bounded. For a fixed $\gamma_j=\gamma\in(\lambda,\lambda_*)$, the squared residuals form a relatively compact set in $C^0(\Omega)$. Weak convergence is uniform on that set, so $J_\gamma(u_j)\to0$ and the fixed-margin stopping condition eventually holds.

For a certified rule with $\gamma_j=\lambda+1/j$, take one node per cell of a measurable partition, weighted by its $m$-measure, and let the maximum cell diameter be $\delta_j\le 2^{-j}$. Given rigorous bounds $S_j$ and $B_j$ on $\|R_j(u_j)\|_\infty$ and its H\"older seminorm, the penalty satisfies
\begin{equation}\label{eq:certified_quadrature}
 J_{\gamma_j}(u_j)\le J_j(u_j)+m(\Omega)S_j B_j\delta_j^\theta .
\end{equation}
Indeed, the seminorm of $R_j(u_j)^2/2$ is at most $S_jB_j$. Taking $S_j,B_j$ from the Sobolev norm and the data keeps them bounded. The right-hand side of~\eqref{eq:certified_quadrature} decays exponentially, and comparison with the polynomially decreasing threshold proves finite termination. These iteration guarantees use exact gradient evaluations for the stated objectives. For a rounded candidate one can test a rigorous upper bound in~\eqref{eq:certified_quadrature} or certify $D_u$ independently, as in the angular-annulus certificate below.

\begin{example}
\label{Sec:Annulus_flow_example}

On the Euclidean annulus $A=\{1\le r\le2\}$, take
$f=-r^2(r-3/2)\partial_r+r^{-1}\partial_\theta$.
This nonvanishing field points inward: $f^r(1)=1/2$ and $f^r(2)=-2$. Its divergence is $b=-4r^2+9r/2$, positive at $r=1$, negative at $r=2$, and zero only at $9/8$.
This annulus illustrates uniform weighted-volume contraction, a distinct objective from convergence of individual trajectories or contraction of their separation~\cite{Rantzer01,Forni14}.

For $\alpha=\gamma=1$ and $\Phi(s)=(s+1)^2/2$, the logarithmic flow is $\psi_t=\diver_g[(b+f(\psi)+1)f]$. At $\psi=0$, its formal initial velocity is $G_0=f(b)+(b+1)b=24r^4-105r^3/2+23r^2+9r/2$. Thus $G_0(1)=-1$, $G_0(2)=65$, and $G_0'(1)=-11$, and $G_0'(2)=469/2$. Whenever compatible boundary data admit a smooth expansion $\psi=tG_0+O(t^2)$, the correction $f^r\psi_r$ is negative on compact subintervals near either boundary, by the signs of the radial velocity and of $G_0'$ there. These are interior first variations, not an existence or convergence proof. On $A$, energy differentiation also produces the boundary term $\int_{\partial A}\Phi'(D_wf)\langle f,\nu\rangle_g\psi_t\,dS_g$.

An independent static certificate is $w_C=\exp[-C(r-1)(2-r)]$, $C\ge0$, for which
\begin{equation}\label{eq:annulus_static_family}
 D_{w_C}f=-4r^2+\frac92r-2Cr^2(r-3/2)^2.
\end{equation}
For $C=13/4$ and $z=r-1\in[0,1]$, the identity $D_{w_C}f+9/8=-z(2-7z+52z^2+52z^3)/8$ and the bound $2-7z+52z^2\ge367/208$ prove $D_{w_C}f\le-9/8$, with equality only at $r=1$. Freezing $w^*=w_{13/4}$ therefore gives $\vol_{(w^*)^{2/n}g}(\varphi_t E)\le e^{-9t/8}\vol_{(w^*)^{2/n}g}(E)$ for all $E\subset A$ and $t\ge0$. This weight is not identified with a flow limit. Its two-sided penalty is nonzero because $D_{w^*}f$ is not identically $-1$.
Figures~\ref{Ex1_fig1}--\ref{Ex1_fig2} show the field, weights, and divergence for static choices of $C$. Maxima in~\eqref{eq:annulus_static_family} are checked at endpoints and interior stationary points.
\begin{figure}[tbp]\centering
\manuscriptfigure{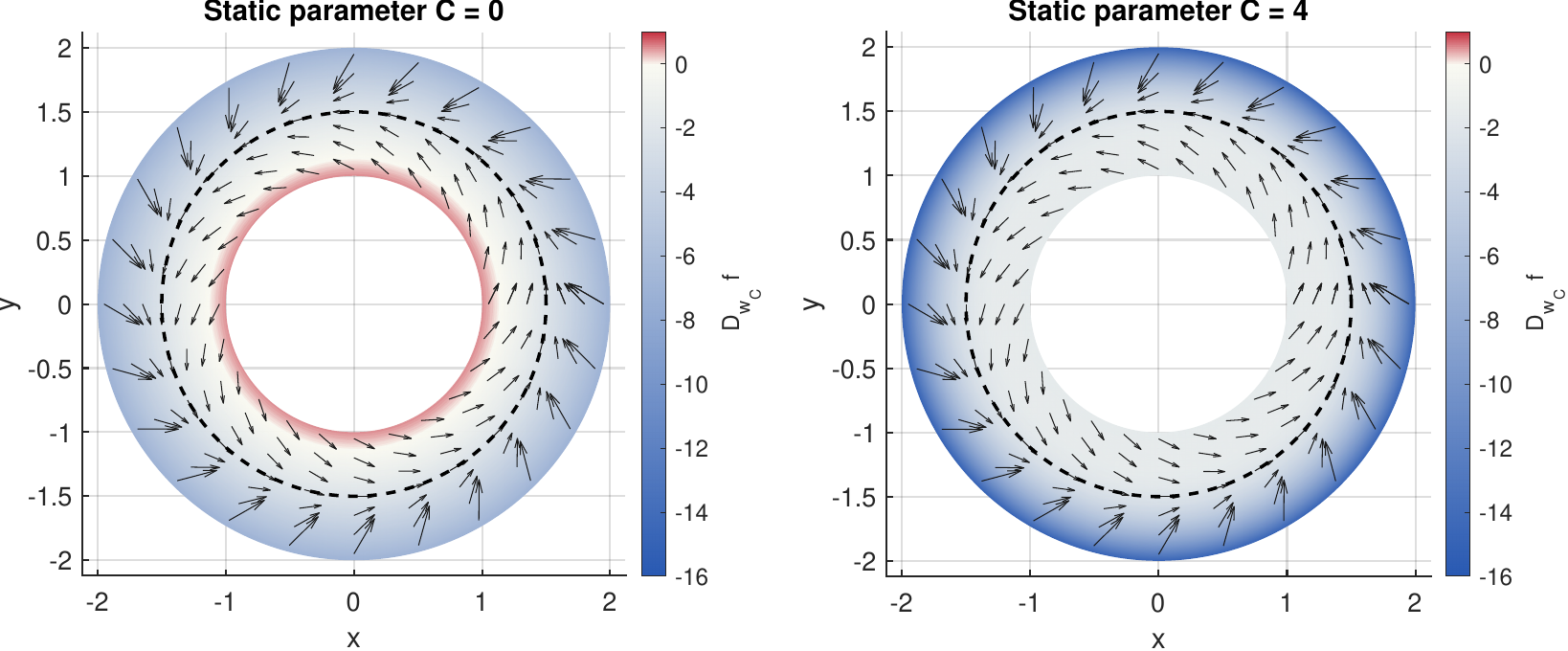}
\caption{Annulus field and exact divergence for $C=0$ and $C=4$. Color denotes $D_{w_C}f$, the dashed circle $r=3/2$ is invariant. Density selection leaves $f$ fixed.}\label{Ex1_fig1}
\end{figure}
\begin{figure}[tbp]\centering
\manuscriptfigure{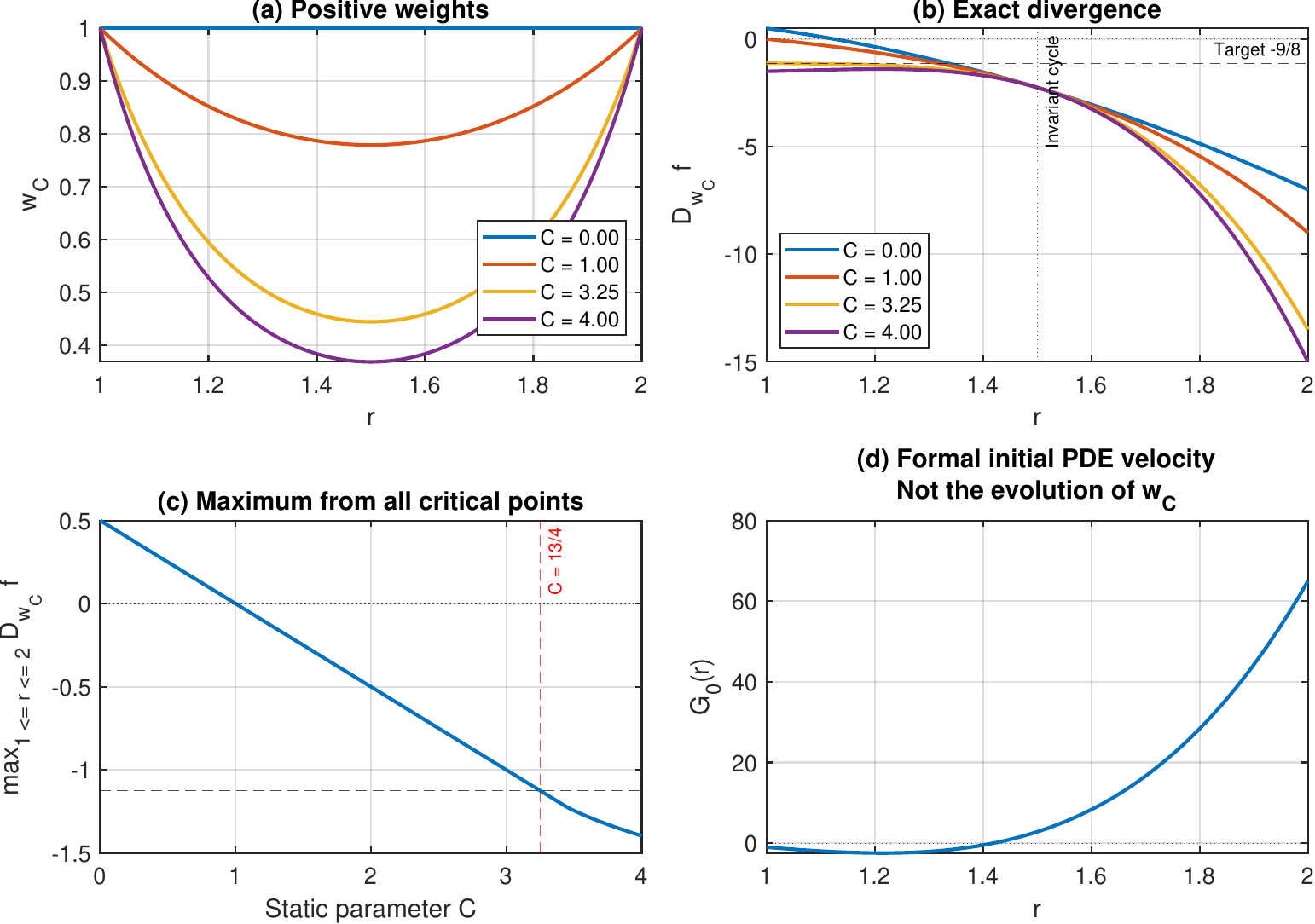}
\caption{Static densities and divergence for $C=0,1,13/4,4$. The exact maximum at $C=13/4$ is $-9/8$. The lower-right panel is the formal initial velocity of the density equation.}\label{Ex1_fig2}
\end{figure}

Set $a=3/2$. The smooth potential $v=-4\log r+a/r$ satisfies $f(v)=(r-a)(4r+a)$ and $D_v=-a^2=-9/4$. The circle $r=a$ is a periodic orbit on which $b=-9/4$. Its invariant probability measure rules out any larger rate by Proposition~\ref{prop:static_feasibility}. Thus $\lambda_*=9/4$, attained by $e^v$. With $\gamma=2$ and target $\lambda=3/2$, the competitor $v$ has zero one-sided penalty and a nonzero two-sided quadratic penalty. Since $H^4(A)\hookrightarrow C^{1,1/2}$, Theorem~\ref{t:main2} gives attainment of every positive rate below $2$ for this choice of $\gamma$. Its stopping test uses $\eta=1/2$ for the target $3/2$. The minimum-norm feasible potential may differ from $v$.

For illustration, restrict to $u=cv_0$, where $v_0=v-v(a)$, and set $M=\|v_0\|_H^2$ using the full norm $\|u\|_H^2=\sum_{j=0}^4\int_A|\nabla^ju|^2dx$, with all ordered Euclidean tensor indices included. For $d=f(v_0)=-9/4-b$, the restricted gradient equation is
\begin{equation}\label{eq:example_annulus_ode}
 M\frac{dc}{d\tau}=-2\pi\int_1^2r\,d(r)(b(r)+cd(r)+2)_+dr
 -\varepsilon Mc,\qquad c(0)=0.
\end{equation}
The computation of~\eqref{eq:example_annulus_ode} in Fig.~\ref{fig:search_computed} uses $\varepsilon M=10^{-2}$ and $\sigma=\tau/M$. The interval $[0,1]$ is invariant. The integral contribution is nonnegative for $0\le c\le1$ because a positive residual requires $b>-9/4$ and hence $d<0$. At $c=1$ the right-hand side is $-\varepsilon M$. At $0\le c\le8/11$ and $1\le r\le21/20$, one has $b>3/10$ and the residual is at least $(3/11)(3/10)+4/11>2/5$. On this radial interval, $-d>5/2$, so the contribution to the right-hand side exceeds $2\pi/20$, greater than $\varepsilon M=10^{-2}$. Thus $c$ reaches $8/11$ in finite time. On $[0,1]$, $D_{cv_0}=(1-c)b-9c/4$ and $\max_A D_{cv_0}=1/2-11c/4$. Thus $c\ge8/11$ certifies $D\le-3/2$ on the entire annulus. This restricted computation does not represent minimization over all of $H^4$.
After freezing the computed density, Fig.~\ref{fig:search_transport} transports $1.02\le r\le1.98$. Its boundaries satisfy $\dot r=-r^2(r-a)$, and its weighted volume $2\pi\int_{r_-(t)}^{r_+(t)}r e^{cv_0(r)}dr$ is bounded by its initial value times $e^{-3t/2}$.
\begin{figure}[tbp]\centering
\manuscriptfigure{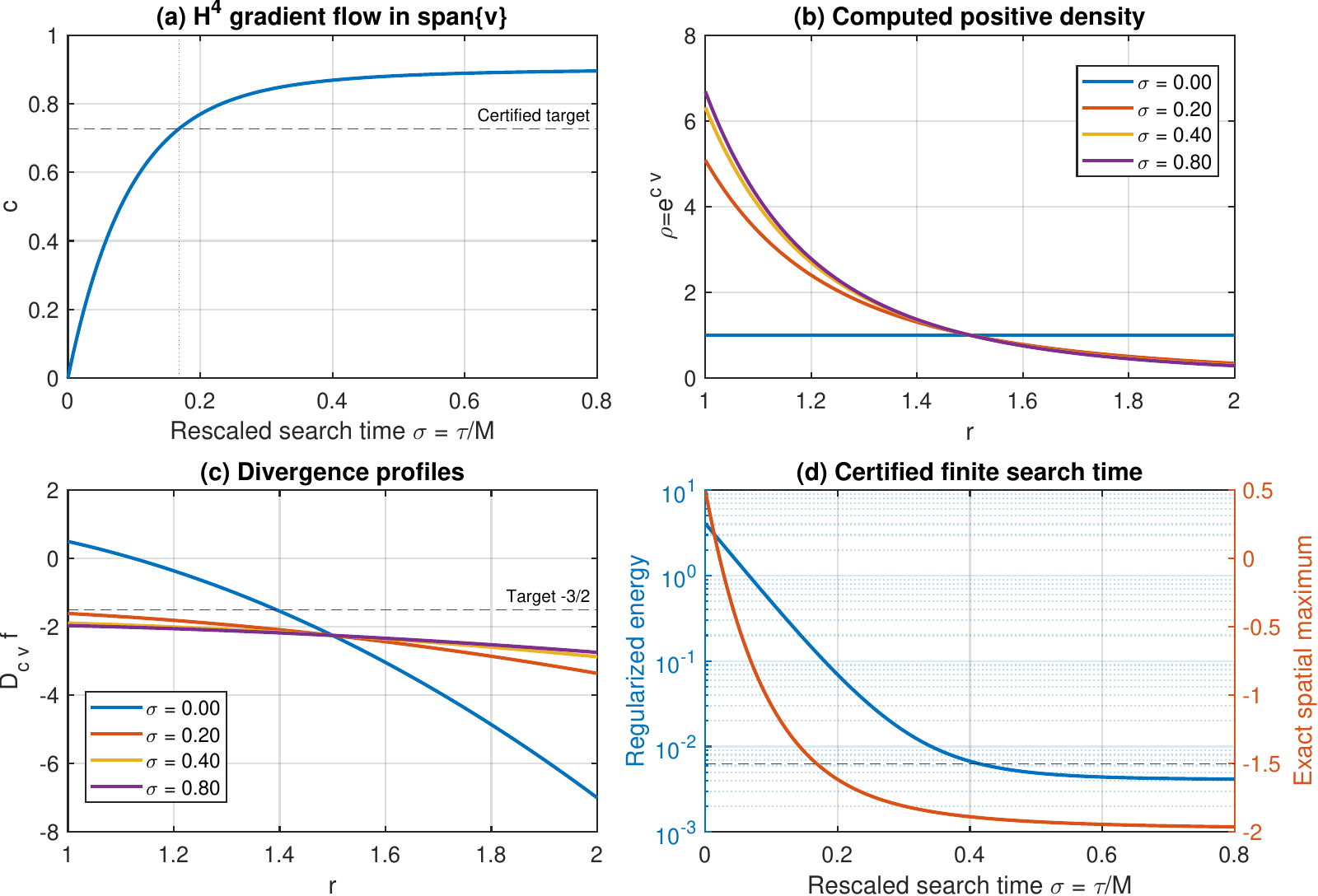}
\caption{Restricted Sobolev search: coefficient, density, divergence, and energy versus $\sigma=\tau/M$. The exact maximum formula certifies rate $3/2$ once $c\ge8/11$.}\label{fig:search_computed}
\end{figure}
\begin{figure}[tbp]\centering
\manuscriptfigure{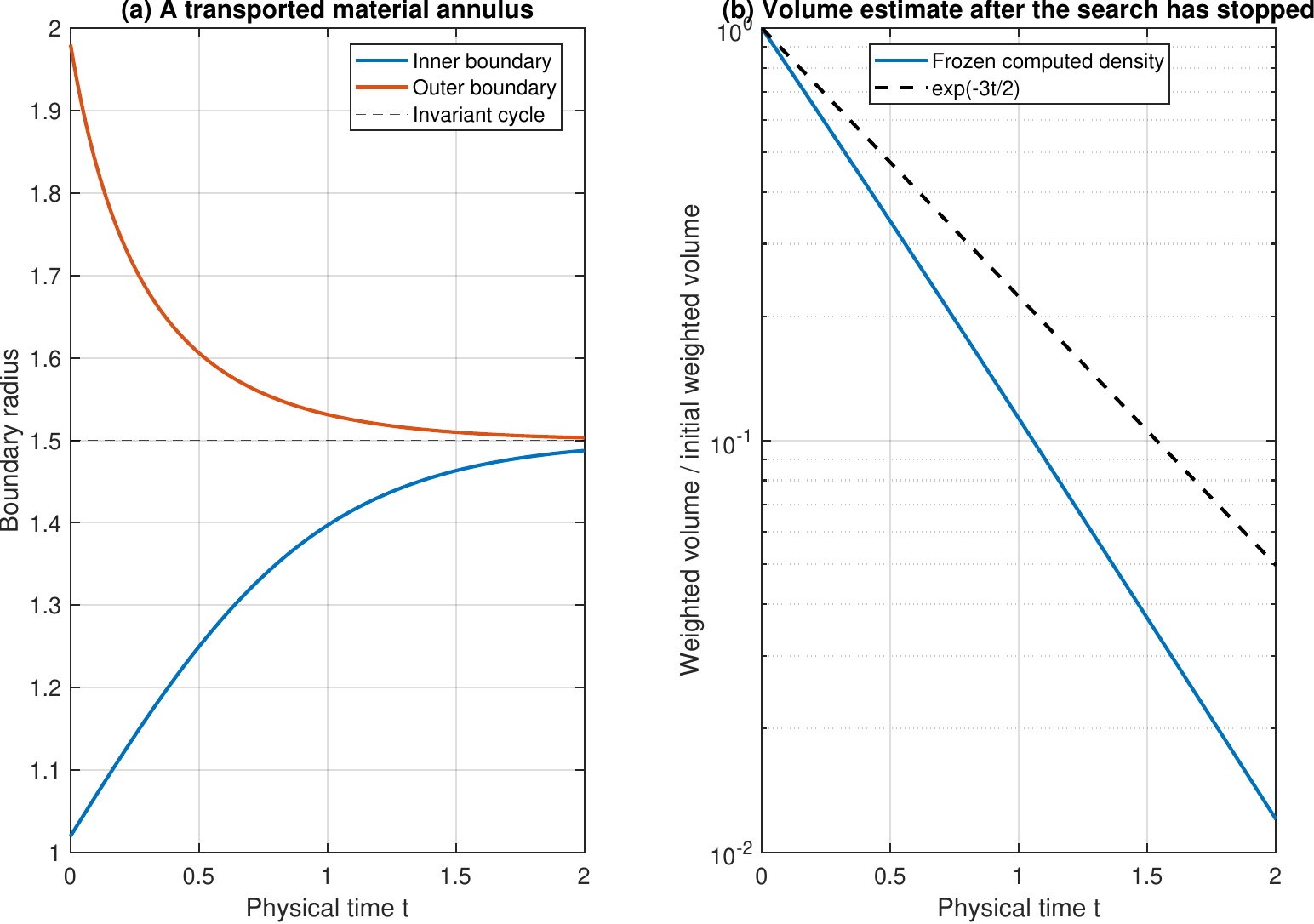}
\caption{Physical transport after freezing the density. Weighted volume is compared with the bound $e^{-3t/2}$, the density remains fixed.}\label{fig:search_transport}
\end{figure}
\end{example}

Keep $f^r=-r^2(r-3/2)$ and take $f^\theta=(1+\tfrac34\sin\theta)/r$. The potential $v(r)-\log(1+\tfrac34\sin\theta)$, with $v$ from Example~\ref{Sec:Annulus_flow_example}, has $D_u\equiv-9/4$, which also identifies $\lambda_*=9/4$ by Proposition~\ref{prop:static_feasibility}. For every radial potential the maximum of $D_u$ on $r=3/2$ is $-7/4$. Thus the target $\lambda=2$ requires angular modes.
The search starts at zero with $\gamma=17/8$ and $\varepsilon=10^{-7}$, using $(2r-3)^i$ times $1,\cos(k\theta),\sin(k\theta)$ for $0\le i\le p$, $1\le k\le m$. The coordinate $H^3$ norm uses $dr\,d\theta/(2\pi)$ and is equivalent to the $H^3$ norm on the annulus. Gradient steps use 32 radial Gauss nodes, 96 equispaced angles and step size $h$ equal to $0.99$ divided by the Lipschitz bound for the quadrature objective. In each rotated quadrant, set $\cos\theta=(1-t^2)/(1+t^2)$ and $\sin\theta=2t/(1+t^2)$ for $0\le t\le1$. The polynomial $r(1+t^2)^{m+1}(D_u-d)$ has nonpositive Bernstein coefficients on subdivided rectangles in the certificate test, which yields $D_u\le d$ everywhere. With checks every 100 steps, both spaces $(p,m)=(3,3)$ and $(4,4)$ first certify $d=-2.01$ at step 800. Their dimensions are 28 and 45 and their search times $800h$ are $112.302$ and $112.285$, respectively. The certificate for the potential with coefficients rounded to nine decimal places is checked in exact rational arithmetic.

On a closed disk of radius less than one, take $f(x,y)=(-x^3,(x-1)y)$. The derivative of $(x^2+y^2)/2$ is $-x^4+(x-1)y^2$, so the field points inward and invariance forces every invariant probability measure to be supported at the origin. Since $b=-1+x-3x^2$, one has $\lambda_*=1$. A $C^1$ potential with $D_u\le-1$ would satisfy $u_x(x,0)\ge x^{-2}-3x^{-1}$ for small $x>0$, which is impossible. For $\delta>0$, the smooth potential $u_\delta=\delta^{-1}\arctan(x/\delta)$ gives $D_{u_\delta}=-1-3x^2+x\delta^2/(x^2+\delta^2)\le-1+\delta/2$. Thus the rates approach the bound even though the bound itself has no $C^1$ certificate.

%%%%%%%%%%%%%%%%%%%%%%%%%%%%%%%%%%%%%%%%%%%%%%%%%%%

\section{Local D-Liouville Flow}
\label{sec:convergence}
\label{sec:flows}

Let $(M,g)$ be a compact oriented Riemannian manifold without boundary, let $f\in\mathfrak{X}(M)$ be a smooth vector field, let $\alpha\in\R$, and let $\Phi:\R\to\R_+$ be smooth and convex. Consider
$\mathcal{J}_{\rm dens}(w)=\int_M \Phi\bigl(D_{w^\alpha}f\bigr)\,dV_g$,
on the space of positive functions $w\in C^\infty(M)$, where $D_{w^\alpha}f$ is defined in~\eqref{eq:weighted_div}. The $L^2$-gradient flow in the variable $\psi=\ln w$, with respect to the standard $L^2$ metric on functions, yields the following evolution equation for $w$:
\begin{equation}
\label{eq:dens_flow}
\partial_t w = \alpha\,w\,\diver_g\!\Bigl(\Phi'\bigl(D_{w^\alpha}f\bigr)\,f\Bigr).
\end{equation}
This equation is called the D-Liouville flow. In the search equations, $t$ denotes geometric evolution time, independent of the physical time of the system. When both processes are considered together, search time is denoted by $\tau$.

For $\psi=\ln w$, $D=\diver_gf+\alpha f(\psi)$. Integration by parts gives
$D\mathcal J_{\rm dens}(\psi)[\eta]=\alpha\int_M\Phi'(D)f(\eta)dV_g=-\alpha\int_M\eta\,\diver_g(\Phi'(D)f)dV_g$.
Thus $\psi_t=\alpha\diver_g(\Phi'(D)f)$, and multiplication by $w$ gives~\eqref{eq:dens_flow}.
Writing the equation directly for $\psi$ gives
\begin{equation}\label{eq:dens_psi}
 \partial_t\psi=\alpha\diver_g\bigl(\Phi'(D_{w^\alpha}f)f\bigr).
\end{equation}
For any $C^2$ penalty and along a smooth solution, the energy identity is
$ \frac{d}{dt}\mathcal J_{\rm dens}(w(t))=-\|\partial_t\psi\|_2^2
 =-\alpha^2\|\diver_g(\Phi'(D_{w^\alpha}f)f)\|_2^2$.
The principal symbol is $\alpha^2\Phi''(D_{w^\alpha}f)(f\cdot\xi)^2$. It has rank at most one and is positive in the $f$ direction where $\alpha\ne0$, $\Phi''>0$, and $f\ne0$.

The one-sided penalty $\Phi(s)=(s+\gamma)_+^2/2$ is $C^{1,1}$, not $C^2$. Classical theorems require a smooth convex approximation or a separate weak formulation. For $\alpha\ne0$, density stationarity means $\diver_g(\Phi'(D)f)=0$, which alone gives no sign constraint. For $\alpha=0$, the flow is stationary for every weight. Nor can the sign of $w_t$ be inferred from $D>0$: derivatives of $D$ and $f$ enter the equation. Quadratic trial weights such as $x^\top Px$ need not form an invariant family. Their evolution requires projection or additional assumptions, as in~\cite{Furtat26_Arc}. Positive $\alpha$ makes $f(\ln w)<0$ decrease divergence, negative $\alpha$ reverses this effect.

\begin{proposition}
\label{th:local_dens}
Let $M$ be a compact Riemannian manifold without boundary, let $f$ be smooth, with zeros allowed, and let $\alpha\in\mathbb R$, $\Phi(s)=\frac12(s+\gamma)^2$, and $\gamma\in\mathbb R$. For every smooth initial datum $\psi_0$, equation~\eqref{eq:dens_psi} has a unique global smooth solution. For initial data $\psi_0\in L^2(M)$, there exists a unique global weak solution in $L^2$. In the general nonlinear case $\Phi''>0$, the principal symbol has rank at most one. The present result does not establish local classical solvability for arbitrary $\Phi$ and arbitrary initial data in $C^{2,\beta}$.
\end{proposition}
\begin{proof}
For the quadratic penalty, write $b=\diver_gf$, $A=\alpha^2\diver_g(f f(\cdot))=\alpha^2(f^2+bf)$, and $h=\alpha\diver_g((b+\gamma)f)$. The equation is $\psi_t=A\psi+h$. Since $\langle u,Au\rangle_2=-\alpha^2\|f(u)\|_2^2$, the closed quadratic form defines a Friedrichs realization generating a contraction semigroup $P_t$. Duhamel's formula $\psi(t)=P_t\psi_0+\int_0^tP_{t-s}h\,ds$ gives the unique global $L^2$ weak solution.

Smoothness cannot be inferred from ellipticity. Instead, the Stratonovich equation $dX_t=\sqrt2\alpha f(X_t)\circ dW_t+\alpha^2b(X_t)f(X_t)dt$, $X_0=x$, has backward generator $A$. Compactness and smooth coefficients give a global smooth stochastic flow with finite derivative moments of every fixed order on bounded time intervals~\cite{Kunita1990}. Thus $\psi(t,x)=\mathbb E\psi_0(X_t^x)+\int_0^t\mathbb E h(X_s^x)ds$ is spatially smooth for smooth data. Uniqueness identifies this with the closed-form semigroup solution, the equation gives all time derivatives. For general $\Phi$, the symbol $\alpha^2\Phi''(D)(f\cdot\xi)^2$ has rank at most one, so this linear argument gives no general nonlinear classical existence theorem.
\end{proof}

\begin{proposition}
\label{th:density_convergence_reg}
Let $M$ be a compact Riemannian manifold with smooth boundary, and let $f\in C^\infty(TM)$, $\alpha\ne0$, and $\gamma>0$. Let $\Phi\in C^2(\mathbb R)$ be convex and nonnegative, with $\Phi(s)=0$ for $s\le-\gamma$, $\Phi(s)>0$ for $s>-\gamma$, and bounded $\Phi''$. For $\varepsilon>0$, set
$\mathcal J_{\rm dens}^\varepsilon(\psi) =\int_M\Phi(\diver_gf+\alpha f(\psi))\,dV_g +\varepsilon\int_M(|\nabla_g\psi|^2+\psi^2)\,dV_g$.
Then $\mathcal J_{\rm dens}^\varepsilon$ has a unique minimizer $\psi_\varepsilon$ on $H^1(M)$. Its $L^2$-gradient flow, understood through the subdifferential of a closed convex functional on $L^2(M)$, exists globally for every $\psi_0\in H^1(M)$ and, at smooth points, has the form
\begin{equation}\label{eq:reg_flow_density}
\partial_t\psi=\alpha\diver_g(\Phi'(D_{e^{\alpha\psi}}f)f)
+2\varepsilon\Delta_g\psi-2\varepsilon\psi,
~~ D_{e^{\alpha\psi}}f=\diver_gf+\alpha f(\psi),
\end{equation}
with the \emph{natural mixed} boundary condition
\begin{equation}\label{eq:reg_natural_boundary}
\alpha\Phi'(D_{e^{\alpha\psi}}f)\langle f,\nu\rangle_g
+2\varepsilon\partial_\nu\psi=0\quad\text{on }\partial M.
\end{equation}
As $t\to\infty$, solutions converge strongly in $H^1(M)$ to $\psi_\varepsilon$. If, additionally, a smooth $\psi^*$ exists such that $\diver_gf+\alpha f(\psi^*)\le-\gamma$ on $M$, then
$0\le\mathcal J_{\rm dens}(\psi_\varepsilon) \le\varepsilon\|\psi^*\|_{H^1}^2\longrightarrow0 $, $(\varepsilon\downarrow0)$.
The last estimate gives smallness of the \emph{integral penalty}. Without further a priori estimates, it does not imply pointwise negativity of the minimizer's divergence.
Set $e_0=\mathcal J_{\rm dens}^\varepsilon(\psi_0)-\mathcal J_{\rm dens}^\varepsilon(\psi_\varepsilon)$. For the corresponding flow and every $t\ge0$, one has $\|\psi(t)-\psi_\varepsilon\|_{L^2}\le e^{-2\varepsilon t}\|\psi_0-\psi_\varepsilon\|_{L^2}$, energy gap $\mathcal J_{\rm dens}^\varepsilon(\psi(t))-\mathcal J_{\rm dens}^\varepsilon(\psi_\varepsilon)\le e_0e^{-4\varepsilon t}$, and $\|\psi(t)-\psi_\varepsilon\|_{H^1}\le\varepsilon^{-1/2}e_0^{1/2}e^{-2\varepsilon t}$.
Thus the regularized flow converges to the unique stationary state from every admissible initial condition.
\end{proposition}
\begin{proof}
The functional is continuous and convex on $H^1$ by boundedness of $\Phi''$, and coercive and uniformly convex because of the regularizer. The direct method gives a unique minimizer. Its extension by $+\infty$ to $L^2\setminus H^1$ is proper and lower semicontinuous, its maximal monotone subdifferential generates the global flow~\cite{Brezis}. Integration by parts gives~\eqref{eq:reg_flow_density} and the full mixed flux~\eqref{eq:reg_natural_boundary}. Convex-flow estimates and uniform convexity give
$\varepsilon\|\psi(t)-\psi_\varepsilon\|_{H^1}^2\le\mathcal J_{\rm dens}^\varepsilon(\psi(t))-\mathcal J_{\rm dens}^\varepsilon(\psi_\varepsilon)\le\|\psi_0-\psi_\varepsilon\|_{L^2}^2/(2t)$.
Comparison with the feasible $\psi^*$ proves the penalty estimate, without pointwise control.
The regularizer has second variation $2\varepsilon\|h\|_{H^1}^2$, so the full functional is $2\varepsilon$-strongly convex in $L^2$. Strong monotonicity yields $\tfrac12\partial_t\|\psi-\psi_\varepsilon\|_2^2\le-2\varepsilon\|\psi-\psi_\varepsilon\|_2^2$. Writing $e=\mathcal J_{\rm dens}^\varepsilon(\psi)-\min\mathcal J_{\rm dens}^\varepsilon$, strong convexity gives $\|\partial^0\mathcal J_{\rm dens}^\varepsilon(\psi)\|_2^2\ge4\varepsilon e$, where $\partial^0$ is the minimal subgradient. The energy identity implies $e'\le-4\varepsilon e$. Finally, $e\ge\varepsilon\|\psi-\psi_\varepsilon\|_{H^1}^2$ gives the last estimate.
\end{proof}

\begin{proposition}
\label{th:nonlinear_density_gap}
Let $M$ be a compact Riemannian manifold, possibly with smooth boundary, let $f$ be a smooth vector field, and let $\alpha\ne0$. Denote by
$H_f^1(M):=\{u\in L^2(M):f(u)\in L^2(M) \text{ in the distributional sense}\}$
the Hilbert space with graph norm $\|u\|_{H_f^1}^2=\|u\|_{L^2}^2+\|f(u)\|_{L^2}^2$. Let $V\subset H_f^1(M)$ be a linear subspace closed in the graph norm $H_f^1(M)$ and dense in $L^2(M)$. Assume that, for some $C_P>0$,
\begin{equation}\label{eq:directional_poincare}
\|u\|_{L^2}^2\le C_P\|f(u)\|_{L^2}^2,
\qquad u\in V.
\end{equation}
The restriction defining $V$ may, for example, impose zero values on the inflow boundary. Fix $\psi_B\in H_f^1(M)$ and the admissible affine space $\psi_B+V$. Let $\Phi\in C^2(\mathbb R)$ be nonnegative and satisfy
$0<m\le\Phi''(s)\le M_0<\infty$ $(s\in\mathbb R)$.
Then the unregularized functional
$E(\psi)=\int_M\Phi(\operatorname{div}_gf+ \alpha f(\psi))\,dV_g$, $\psi\in\psi_B+V$,
has a unique minimizer $\psi_*\in\psi_B+V$. Its negative $L^2$-gradient flow in the subdifferential sense exists for every initial $\psi_0\in\psi_B+V$ and all $t\ge0$, and satisfies
\begin{align}
\|\psi(t)-\psi_*\|_{L^2}
&\le e^{-\kappa t}\|\psi_0-\psi_*\|_{L^2},\label{eq:nonlin_L2_rate}\\
E(\psi(t))-E(\psi_*)
&\le e^{-2\kappa t}[E(\psi_0)-E(\psi_*)],
\label{eq:nonlin_energy_rate}
\end{align}
where $\kappa=m\alpha^2/C_P$. Moreover,
\begin{equation}\label{eq:nonlin_graph_rate}
\|f(\psi(t)-\psi_*)\|_{L^2}^2
\le\frac{2}{m\alpha^2}
[E(\psi_0)-E(\psi_*)]e^{-2\kappa t}.
\end{equation}
In particular, the solution converges exponentially to $\psi_*$ in the graph norm $H_f^1$. The flow satisfies the variational identity
\[
 \langle\partial_t\psi,v\rangle_{L^2}
   =-\alpha\int_M\Phi'(\diver_g f+\alpha f(\psi))\,f(v)\,dV_g,\qquad v\in V,
\]
for almost every positive time. If $C_c^\infty(\operatorname{int}M)\subset V$, this gives~\eqref{eq:dens_psi} in the interior. At a boundary portion where smooth variations in $V$ have arbitrary supported traces, smooth solutions satisfy the natural flux condition $\alpha\Phi'(D_{w^\alpha}f)\langle f,\nu\rangle_g=0$.
\end{proposition}
\begin{proof}
The distributional operator $u\mapsto f(u)$ is closed, and the bounds on $\Phi''$ give continuity and quadratic growth in the graph norm. For $u\in V$, uniform convexity and~\eqref{eq:directional_poincare} give $E(\psi+u)\ge E(\psi)+DE(\psi)[u]+(m\alpha^2/2)\|f(u)\|_2^2\ge E(\psi)+DE(\psi)[u]+(\kappa/2)\|u\|_2^2$. This also gives graph-norm coercivity. The direct method yields a unique minimizer. Bounded sublevels and graph-closedness of $V$ show that extension by $+\infty$ outside $\psi_B+V$ is lower semicontinuous on $L^2$. Its maximal monotone subdifferential therefore generates a global semigroup~\cite{Brezis}.

Strong monotonicity gives $\tfrac12\partial_t\|\psi-\psi_*\|_2^2\le-\kappa\|\psi-\psi_*\|_2^2$. For the minimal subgradient $\xi=-\psi_t$, strong convexity gives $E(\psi)-E(\psi_*)\le\|\xi\|_2^2/(2\kappa)$, while $E'=-\|\xi\|_2^2$. Gronwall proves~\eqref{eq:nonlin_L2_rate}--\eqref{eq:nonlin_energy_rate}. Finally, $DE(\psi_*)[u]=0$ for $u\in V$ implies $E(\psi)-E(\psi_*)\ge(m\alpha^2/2)\|f(\psi-\psi_*)\|_2^2$, giving~\eqref{eq:nonlin_graph_rate}. The first variation gives the variational identity. Compactly supported variations yield the interior equation, and testing the remaining boundary term with arbitrary smooth traces on a free portion gives the flux condition there.
\end{proof}

The assumption on interior variations matters. On $[0,1]$ with $f=\partial_x$, $\alpha=1$, and $\Phi(s)=s^2/2$, take $V=\{v\in H^1:v(1)=0,\ v(0)=\int_0^1v\,dx\}$. This is graph-closed and satisfies the directional Poincar\'e inequality. It is dense in $L^2$ because a bump supported near $0$ can enforce the trace constraint with arbitrarily small $L^2$ norm. With $\psi_B=x-x^2/2$, one has $DE(\psi_B)[v]=-v(0)+\int_0^1v\,dx=0$ on $V$. The smooth stationary state $\psi_B$ has $\psi_B''=-1$, so the variational equation in this case does not give the heat equation in the interior.

\begin{example}
\label{ex:nonlinear_density_interval}
Attainment of strictly negative divergence can be illustrated as follows. Let $M=[0,L]$, $f=\partial_x$, $\alpha\ne0$, $\Phi(s)=\tfrac12(s+\gamma)^2$, $\gamma>0$, and $V=\{u\in H^1(0,L):u(0)=0\}$. The fundamental theorem of calculus gives~\eqref{eq:directional_poincare} with $C_P=L^2/2$. For $\psi_B=0$ on $V$, the exact minimizer is $\psi_*(x)=-\gamma x/\alpha$ and $D_{e^{\alpha\psi_*}}f =\alpha\psi_*'=-\gamma$ on $[0,L].$ Here the flow equation is linear:
$\partial_t\psi=\alpha^2\partial_{xx}\psi$,
$\psi(0,t)=0$,
$\alpha\psi_x(L,t)+\gamma=0$.
Subtracting $\psi_*$ gives the heat equation with a homogeneous Dirichlet condition at the left endpoint and a homogeneous Neumann condition at the right endpoint. Its spectral gap is $\alpha^2\pi^2/(4L^2)$. Thus $\psi(t)$ converges exponentially to $\psi_*$ in $H^1$. For compatible smooth data, convergence in $C^1$ follows from smoothing for every $t\ge t_0>0$. The positive limiting density $e^{\alpha\psi_*}=e^{-\gamma x}$ defines a strictly dissipative volume form. This example demonstrates actual attainment by the flow. However, transport estimates on the interval apply only until trajectories exit through the right endpoint: the field $\partial_x$ does not make the entire interval positively invariant.

On $[0,1]$, let $f=\partial_x$, $\alpha=1$, and $\Phi(D)=(D+1)^2/2$. The mixed-boundary heat flow $u_\tau=u_{xx}$, $u(0,\tau)=0$, $u_x(1,\tau)=-1$, has the exact compatible solution $u=-x+e^{-\pi^2\tau/4}\sin(\pi x/2)$. Its divergence $D_u=-1+(\pi/2)e^{-\pi^2\tau/4}\cos(\pi x/2)$ is initially positive near zero. Figure~\ref{fig:heat_exact} shows that $D_u\le-1/2$ everywhere exactly when $\tau\ge4\log\pi/\pi^2$. A directional Poincar\'e constant is $C_P=1/2$, the sharp mixed-boundary spectral gap is $\pi^2/4$. Physical transport applies only until exit at the right endpoint, unlike the invariant annulus.
\begin{figure}[tbp]\centering
\manuscriptfigure{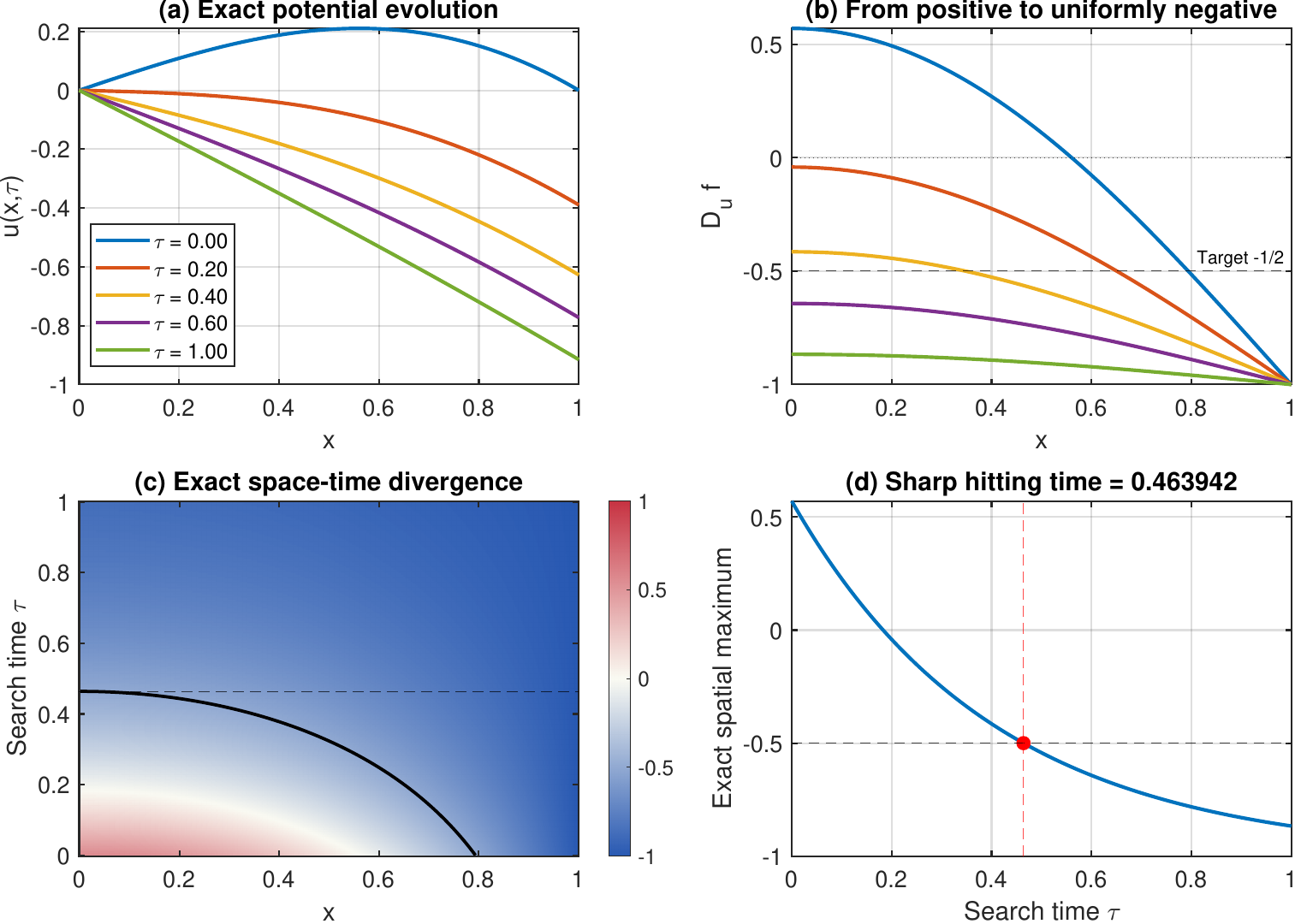}
\caption{Exact heat-flow potential, divergence, space--time map, and maximum. The target $-1/2$ is attained everywhere at $\tau=4\log\pi/\pi^2$. The physical domain is not invariant.}\label{fig:heat_exact}
\end{figure}
\end{example}

The mixed boundary data differ from the natural condition~\eqref{eq:reg_natural_boundary}. For a nonnegative convex penalty with zero set $(-\infty,-\gamma]$, a global unregularized orbit bounded in $H^1$ with penalty tending to zero has a weak subsequential $H^1$ limit of zero penalty and divergence at most $-\gamma$ almost everywhere. Smoothness and integrability of its exponential require further bounds. A positive $C^1$ limiting density on a closed manifold would contradict~\eqref{eq:stokes_obstruction}. Weak convergence alone does not justify applying Stokes, and energy monotonicity alone does not imply a zero minimum. Convergence for a Yamabe functional with $H^1$ regularization or for a metric functional with an added $\varepsilon\int|\Rm|^2dV$ requires a separate argument.

\begin{proposition}
Let $\Omega$ be a compact positively invariant domain with smooth boundary, let $\alpha>0$, and let $\Phi\in C^2(\mathbb R)$ be nonnegative and convex, with zero set $(-\infty,-\gamma]$ for some $\gamma>0$. Suppose the admissible class of smooth positive weights contains $w_0$ satisfying $D_{w_0^\alpha}f\le-\gamma$ on $\Omega$. Then
$\inf_w\mathcal J_{\rm dens}(w)=0$, $\mathcal J_{\rm dens}(w) =\int_\Omega\Phi(D_{w^\alpha}f)\,dV_g$.
Every smooth positive global minimizer $\hat w$ satisfies $D_{\hat w^\alpha}f\le-\gamma$ throughout $\Omega$ and hence gives the volume estimate~\eqref{eq:vol_exp_D} in the conformal metric $\hat w^{2\alpha/n}g$. The same conclusion holds for a positive $C^1$ limit of the density flow if that limit is a global minimizer. A flow on a domain with boundary must additionally have a compatible boundary formulation.
\end{proposition}
\begin{proof}
The feasible competitor has zero energy. At any global minimizer, nonnegativity forces $\Phi(D)=0$ almost everywhere, hence everywhere by continuity, including the boundary. Freezing the density gives~\eqref{eq:vol_exp_D}. Proposition~\ref{prop:static_feasibility} supplies a competitor $w_0=e^{u/\alpha}$ whenever $0<\gamma<\lambda_*(\Omega,f)$ in the unrestricted class. No convergence assertion follows from this comparison.
\end{proof}

%%%%%%%%%%%%%%%%%%%%%%%%%%%%%%%%%%%%%%%%%%%%%%%%%%%

\section{D-Yamabe Flow and the Scaling Obstruction}
\label{sec:yamabe}

For $w=e^{n\phi}$, the metric $\widetilde g_\alpha=e^{2\alpha\phi}g$ has scalar curvature~\cite{Besse}
\[
 \widetilde R_\alpha=e^{-2\alpha\phi}\bigl[R_g-2\alpha(n-1)\Delta_g\phi
             -\alpha^2(n-1)(n-2)|\nabla_g\phi|^2\bigr].
\]
Let $(M,g)$ be a compact oriented Riemannian manifold without boundary of dimension $n\ge3$, let $f$ be smooth, and assume $\alpha\ne0$, $\zeta>0$, and convex $\Phi\in C^\infty(\mathbb R)$. Here $\zeta$ is a penalty coefficient, not a volume-contraction rate. For $w=e^{n\phi}$ and $\widetilde g_\alpha=e^{2\alpha\phi}g$, consider
\begin{equation}\label{eq:F_alpha_corr}
\mathcal F_\alpha(\phi)=\int_M\widetilde R_\alpha\,dV_{\widetilde g_\alpha}
+\zeta\int_M\Phi(D_{w^\alpha}f)\,dV_g.
\end{equation}
Its negative $L^2(dV_g)$-gradient flow is
\begin{equation}\label{eq:YF_exact_corr}
\partial_t\phi=-(n-2)\alpha e^{n\alpha\phi}\widetilde R_\alpha
+\zeta\alpha n\,\diver_g\bigl(\Phi'(D_{w^\alpha}f)f\bigr).
\end{equation}
This is a Yamabe-type equation with a divergence penalty. For $\zeta=0$, it is generally not classical Yamabe flow with its usual time parametrization.

For $\delta\widetilde g=2\alpha\eta\widetilde g$, the total-scalar-curvature variation is $(n-2)\alpha\int_M\widetilde R_\alpha\eta e^{n\alpha\phi}dV_g$. Since $D=\diver_gf+\alpha n f(\phi)$, the penalty variation is $-\zeta\alpha n\int_M\eta\,\diver_g(\Phi'(D)f)dV_g$. Their sum gives~\eqref{eq:YF_exact_corr}, no volume constraint is imposed, so no mean-curvature subtraction appears.
The fixed-metric gradient identity gives
$\partial_t\mathcal F_\alpha(\phi(t))=-\|\nabla_{L^2}\mathcal F_\alpha(\phi(t))\|_2^2$.
If the evolving $L^2(dV_{\widetilde g_\alpha})$ metric on conformal factors is used, the same functional yields
$\partial_t\phi=-(n-2)\alpha\widetilde R_\alpha +\zeta\alpha n e^{-n\alpha\phi} \diver_g(\Phi'(D_{w^\alpha}f)f).$
For $\zeta=0$, this equation induces unnormalized Yamabe flow for $\widetilde g_\alpha$ after a constant time rescaling. In contrast, flow~\eqref{eq:YF_exact_corr} uses the fixed measure $dV_g$ and has a different multiplier on its right-hand side.

\begin{theorem}
\label{prop:yamabe_scaling_obstruction}
Let $M$ be compact and without boundary, let $n\ge3$ and $\alpha\ne0$, and let $\mathcal F_\alpha$ be given by~\eqref{eq:F_alpha_corr}. Every smooth stationary state $\phi_*$ of its negative $L^2(dV_g)$-gradient flow satisfies
$\int_M R_{e^{2\alpha\phi_*}g}\,dV_{e^{2\alpha\phi_*}g}=0$.
Moreover, the following strict alternative holds:

1) either all $\phi_*+c$, $c\in\mathbb R$, are stationary, and the linearization has a zero eigenvalue in the constant direction,

2) or the Hessian of $\mathcal F_\alpha$ at $\phi_*$ is indefinite, and the stationary state is linearly unstable.

In particular, strictly negative spectrum of the linearization, as required for exponential stability of an individual stationary state, is impossible for the original unnormalized flow on the full space of conformal factors.
\end{theorem}
\begin{proof}
Write
$A(\phi)=\int_M R_{e^{2\alpha\phi}g}\,dV_{e^{2\alpha\phi}g}$
and $J(\phi)=\int_M\Phi(\operatorname{div}_gf+
\alpha n f(\phi))\,dV_g$. Adding a constant $c$ gives
$A(\phi+c)=e^{(n-2)\alpha c}A(\phi),\qquad J(\phi+c)=J(\phi).$
Stationarity and variation in the constant direction imply $(n-2)\alpha A(\phi_*)=0$, hence $A(\phi_*)=0$. Therefore, in the constant direction $D^2\mathcal F_\alpha(\phi_*)[1,1]=0$, while for arbitrary $\eta$,
$D^2\mathcal F_\alpha(\phi_*)[1,\eta] =(n-2)\alpha\,DA(\phi_*)[\eta].$
If $DA(\phi_*)\equiv0$, stationarity $D\mathcal F_\alpha(\phi_*)=0$ gives $DJ(\phi_*)=0$. The scaling identities then imply that $\phi_*+c$ is stationary for all $c$, so constants belong to the kernel of the Hessian. If instead $DA(\phi_*)[\eta]\ne0$, write $B=D^2\mathcal F_\alpha(\phi_*)$. Then $B[\eta+t,\eta+t]=B[\eta,\eta]+2t(n-2)\alpha DA(\phi_*)[\eta]$ takes both signs as $t$ varies. At a critical point, the Hessian is a formally self-adjoint second-order elliptic operator relative to the fixed $L^2(dV_g)$ metric. Its self-adjoint realization on compact $M$ has compact resolvent. A direction on which its quadratic form is negative therefore implies a negative eigenvalue of the Hessian, and hence a positive eigenvalue of the linearized negative gradient flow. The stationary state is consequently linearly unstable, proving the alternative.
\end{proof}

To remove the constant direction, restrict the energy to the space
$X_0^{2,\beta}=\{\phi\in C^{2,\beta}(M):\langle\phi\rangle=0\}, \qquad P v=v-\langle v\rangle,$
where $\langle v\rangle=\vol_g(M)^{-1}\int_Mv\,dV_g$, and the normalized flow
\begin{equation}
\label{eq:normalized_Yamabe_general}
 \partial_t\phi=P F_Y(\phi),\qquad \langle\phi(0)\rangle=0,
\end{equation}
where $F_Y$ is the right-hand side of~\eqref{eq:YF_exact_corr}. This is the negative $L^2(dV_g)$ gradient of $\mathcal F_\alpha$ restricted to the mean-zero hyperplane.
Integration by parts gives
$ A(\phi)=\int_M e^{(n-2)\alpha\phi}
   \bigl[R_g+(n-1)(n-2)\alpha^2|\nabla_g\phi|^2\bigr]\,dV_g$.
Since the divergence term in $F_Y$ has zero mean, the spatially constant correction in $PF_Y$ is determined by $A(\phi)$ and involves at most first derivatives.

\begin{proposition}
\label{th:local_YF}
Let $(M,g)$ be compact and without boundary, let $n\ge3$, let $f$ be smooth, and assume $\alpha\ne0$, $\zeta>0$, $\Phi\in C^\infty(\mathbb R)$, and $\Phi''\ge0$. For every $\phi_0\in C^{2,\beta}(M)$, $0<\beta<1$, equation~\eqref{eq:YF_exact_corr} has a unique local classical solution in the standard parabolic class $C^{2+\beta,1+\beta/2}$.
The same conclusion holds for~\eqref{eq:normalized_Yamabe_general} on the mean-zero subspace. In either equation a uniform $C^{2,\beta}$ bound up to a finite maximal time implies continuation.
\end{proposition}
\begin{proof}
Substitution of the conformal scalar-curvature transformation into~\eqref{eq:YF_exact_corr} gives the principal symbol
\begin{equation}\label{eq:main_term}
\sigma(\xi)=2(n-2)(n-1)\alpha^2e^{(n-2)\alpha\phi}|\xi|_g^2
+\zeta\alpha^2n^2\Phi''(D_{w^\alpha}f)(f\cdot\xi)^2.
\end{equation}
On a compact manifold with initial data $\phi_0\in C^{2,\beta}$, the first term is uniformly positive for a short time, while the second is nonnegative. The coefficients depend smoothly on $\phi$ and $\nabla\phi$, and the remaining terms have order at most one. Standard local theory for uniformly parabolic quasilinear equations gives the result.
The projection adds a smooth operator of order at most one, so the same local quasilinear theory applies, and the mean is preserved by the equation. Under the stated bound the time derivative is uniformly bounded in $C^0$, yielding a unique $C^0$ limit at the putative maximal time. Compactness upgrades this to convergence in $C^{2,\beta'}$ for $0<\beta'<\beta$. Uniform parabolicity and interior estimates on a terminal time interval give uniform higher spatial bounds, so the solution restarts in the original regularity class.
\end{proof}

\begin{proposition}
\label{th:normalized_yamabe_torus}
Let $M=\mathbb T^n$, $n\ge3$, carry the flat metric $g_0$. Fix $\alpha\ne0$, $\zeta>0$, $\Phi(s)=s^2/2$, and a smooth vector field $F$ equal to $-\sum_i x^i\partial_{x^i}$ on a coordinate ball $U$ centered at $x_*$. Set $f_\delta=\delta F$. For sufficiently small $\delta>0$, the normalized equation~\eqref{eq:normalized_Yamabe_general} has a unique smooth mean-zero stationary factor $\phi_\delta$ in a fixed neighborhood of zero, with $\phi_\delta=O(\delta^2)$ in $C^{2+\beta}$. All sufficiently close mean-zero data in $h^{2+\beta}$ generate global solutions converging exponentially to $\phi_\delta$. On a smaller positively invariant ball $U_*\Subset U$, the metric $g_\delta=e^{2\alpha\phi_\delta}g_0$ has $\diver_{g_\delta}f_\delta<0$.
\end{proposition}
\begin{proof}
On the mean-zero subspaces of $h^{2+\beta}(M)$ and $h^\beta(M)$, consider $G(\phi,\delta)=P F_Y(\phi,f_\delta)$. One has $G(0,0)=0$, and its derivative with respect to $\phi$ at $(0,0)$ is $D_\phi G(0,0)=2(n-1)(n-2)\alpha^2\Delta_{g_0}.$ The Laplacian is invertible on mean-zero functions and has strictly negative spectrum. The implicit function theorem therefore yields a smooth branch of stationary factors $\phi_\delta$. Since the penalty term at $\phi=0$ with $f_\delta=\delta F$ is of order $O(\delta^2)$, one obtains $\phi_\delta=O(\delta^2)$. The projection $P$ preserves the principal part, as shown above. On the mean-zero subspace, the linearization remains sectorial, and its strict spectral gap persists for small $\delta$. The principle of linearized stability for quasilinear parabolic equations~\cite{PrussSimonett2016} gives global existence and exponential convergence for nearby initial data. Finally, $\operatorname{div}_{g_\delta}f_\delta
=\delta\operatorname{div}_{g_0}F+
\alpha n\delta F(\phi_\delta)$, whose last term is of order $O(\delta^3)$. On a sufficiently small ball around $x_*$, the field $F=-x$ has strictly negative divergence and points inward along the boundary. Taking $\delta>0$ sufficiently small yields strictly negative divergence in $g_\delta$ on the same smaller ball.
\end{proof}

\begin{proposition}
\label{th:global_normalized_Yamabe}
Let $M$ be compact and without boundary, let $n\ge3$, $\alpha\ne0$, and $\zeta>0$, and let $\Phi$ be real analytic and convex with $\Phi''\ge0$. Let $\phi(t)$ be the maximal classical solution of~\eqref{eq:normalized_Yamabe_general} with $\phi_0\in X_0^{2,\beta}$ and $0<\beta<1$. Assume that
\begin{equation}
\label{eq:Yamabe_global_apriori}
 \sup_{0\le t<T_{\max}}\|\phi(t)\|_{C^{2,\beta}(M)}<\infty.
\end{equation}
Then $T_{\max}=\infty$, and there exists a smooth stationary solution $\phi_\infty$ of the normalized flow such that
$\phi(t)\longrightarrow\phi_\infty \qquad\text{in }C^k(M)\quad\text{for every }k\ge0.$
For some \L{}ojasiewicz--Simon exponent $\theta\in(0,1/2]$, the estimates of Lemma~\ref{th:abstract_LS} hold: if $0<\theta<1/2$, convergence in $L^2$ has the algebraic rate stated there, while for $\theta=1/2$ it is exponential. In particular, if the Hessian of the restriction $\mathcal F_\alpha|_{\langle\phi\rangle=0}$ at $\phi_\infty$ is nondegenerate, convergence in $L^2$ is exponential. Parabolic smoothing transfers the exponential rate to every fixed stronger norm after positive time, possibly with a smaller exponent. The same conclusion holds in the Morse--Bott case with optimal exponent $\theta=1/2$.
\end{proposition}

\begin{proof}
Proposition~\ref{th:local_YF} gives global existence under~\eqref{eq:Yamabe_global_apriori}. Schauder estimates on translated cylinders then give uniform $C^m$ bounds on every positive-time tail, so the tail is precompact in $X_0^{2,\beta_0}$ for $0<\beta_0<\beta$. The restricted gradient identity is
\begin{equation}\label{eq:normalized_Yamabe_energy_identity}
 \partial_t\mathcal F_\alpha(\phi(t))
 =-\|P\nabla_{L^2}\mathcal F_\alpha(\phi(t))\|_2^2.
\end{equation}
The energy and projected gradient are real analytic on the H\"older chart with exponent $\beta_0$. The derivative of the projected gradient acts from $C_0^{2,\beta_0}$ to $C_0^{0,\beta_0}$ and extends to a bounded operator from $H_0^2$ to $L_0^2$. The extension depends continuously on the factor in operator norm, since its coefficients depend continuously on the factor in $C^{2,\beta_0}$ and the projection is bounded. At a stationary factor the Hessian is formally self-adjoint and elliptic, with positive principal symbol~\eqref{eq:main_term}. The projection is a finite-rank correction. Elliptic regularity identifies the kernel and the orthogonality conditions for the range in the two mean-zero realizations, which are both Fredholm of index zero. The Hilbert-norm version of the \L{}ojasiewicz--Simon inequality (Theorem~3 of~\cite{FeehanMaridakis2017}) now applies with the $L^2$ gradient norm. Lemma~\ref{th:abstract_LS}, with $X=X_0^{2,\beta_0}$ and $H=L_0^2(M,dV_g)$, gives convergence to one critical factor and the stated rates. Higher bounds and interpolation give $C^k$ convergence. For a nondegenerate critical point the Lyapunov--Schmidt reduction has a zero-dimensional kernel. In the Morse--Bott case the reduced energy is locally constant because the critical manifold projects locally onto the kernel. In either case the quadratic estimate transverse to the kernel gives $\theta=1/2$ in the $L^2$ inequality~\cite{FeehanMaridakis2017}. Parabolic estimates transfer exponential decay to fixed stronger norms.
\end{proof}

A nonzero real-analytic penalty cannot vanish on the whole half-line $s\le-\gamma$. The convergence argument above therefore concerns analytic penalties. The squared positive-part penalty in Section~\ref{s:search2} has its own convex-flow proof.

In a positively invariant class, a $C^{2,\beta}$ bound on the initial energy sublevel yields~\eqref{eq:Yamabe_global_apriori} by~\eqref{eq:normalized_Yamabe_energy_identity}, so Proposition~\ref{th:global_normalized_Yamabe} applies. Boundedness of every sublevel gives convergence for every smooth mean-zero initial factor in that class. On the full mean-zero space, take $\phi_k=k^{-3/2}\chi\sin(kx_1)$ in a chart with nonzero smooth cutoff $\chi$ and subtract its mean. Then $\phi_k\to0$ in $C^1$, so $\mathcal F_\alpha(\phi_k)\to\mathcal F_\alpha(0)$ by the first-derivative formula for $A$, but $\|\phi_k\|_{C^2}\to\infty$. Thus every sublevel above $\mathcal F_\alpha(0)$ is unbounded in $C^{2,\beta}$.

With $f=0$ on a round sphere with scalar curvature $R_0>0$, the spatially constant solution of~\eqref{eq:YF_exact_corr} from $\phi(0)=0$ satisfies $\dot\phi=-cR_0e^{c\phi}$, where $c=(n-2)\alpha$. Hence $\phi(t)=-c^{-1}\log(1+c^2R_0t)$ and $\widetilde g(t)=(1+c^2R_0t)^{-2/(n-2)}g_0$. The search is global and the metric degenerates as $t\to\infty$.

The mean-zero projection fixes constant factors on a connected round sphere in this example. On a disconnected manifold the normalized equation itself can degenerate. Take a round sphere and a flat $n$-torus of equal background volumes, with $f=0$. Factors $a(t)$ and $-a(t)$ on the two components solve the normalized equation when $a'=-(cR_0/2)e^{ca}$ and $a(0)=0$. Thus $e^{-ca(t)}=1+c^2R_0t/2$. The spherical metric shrinks and the torus metric expands, so this smooth global search has no Riemannian limit.

\begin{example}
\label{Sec:Torus_example}
Consider conformal deformation without a stationary limit.
On $\mathbb T^3=\mathbb R^3/(2\pi\mathbb Z)^3$ with $g_0=dx^2+dy^2+dz^2$, let $f=(1+\sin x)\partial_x+\partial_y+\partial_z$. Then $b=\cos x$, while classical Yamabe flow leaves $g_0$ fixed. For the D-Yamabe flow with $\Phi(s)=(s+\gamma)^2/2$, $\gamma>0$, $\alpha,\zeta>0$, and $\phi(0)=0$, Proposition~\ref{th:local_YF} gives local existence. Symmetry reduces the equation to the $x$ circle, with $D=\cos x+3\alpha(1+\sin x)\phi_x$. Its initial velocity is $\phi_t|_0=3\zeta\alpha[-\sin x(1+\sin x)+\cos^2x+\gamma\cos x]$. At $x=0$, this equals $3\zeta\alpha(1+\gamma)>0$, whereas $\phi_x(0,t)=-3\zeta\alpha t+O(t^2)$ and $D(0,t)=1-9\zeta\alpha^2t+O(t^2)$. Divergence therefore decreases locally for small positive time, without a finite-time negativity conclusion, see Fig.~\ref{Ex2_fig1}.

No smooth stationary factor exists. Integrating the stationary equation gives zero total scalar curvature, but the conformal formula on the flat torus gives $\int\widehat R_\alpha\,dV_{\widehat g_\alpha}=2\alpha^2\int e^{\alpha\hat\phi}|\nabla\hat\phi|^2dV_{g_0}$. Hence $\hat\phi$ is constant. Stationarity would then require $\partial_x[(\cos x+\gamma)(1+\sin x)]=0$, whose value at zero is $1+\gamma$. Thus no global smooth trajectory can converge in $C^2$ to a smooth stationary metric. Stokes' identity also excludes negative divergence on the whole torus. Density selection has the same static possibilities for divergence, although it leaves curvature fixed.

For $\alpha=\zeta=\gamma=1$, put $a=1+\sin x$ and $D=\cos x+3a\phi_x$. Figure~\ref{Ex2_fig2} integrates the nonlinear equation~\eqref{eq:example_yamabe_reduced},
\begin{equation}\label{eq:example_yamabe_reduced}
 \phi_\tau=e^\phi(4\phi_{xx}+2\phi_x^2)+3\partial_x[(D+1)a],
 \qquad \phi(x,0)=0.
\end{equation}
The energy per unit $(y,z)$ coordinate area is $\mathcal E=2\int e^\phi\phi_x^2dx+\tfrac12\int(D+1)^2dx$, and $\mathcal E'=-\int\phi_\tau^2dx$. Fourier collocation uses $64$ and $128$ equispaced nodes on $[0,2\pi)$, with products evaluated at the nodes. 
Time integration on $[0,0.2]$ was performed using \texttt{ode15s}, with relative and absolute tolerances of $10^{-10}$ and $10^{-12}$, respectively, and a maximum time step of $0.005$.

\begin{figure}[tbp]\centering
\manuscriptfigure{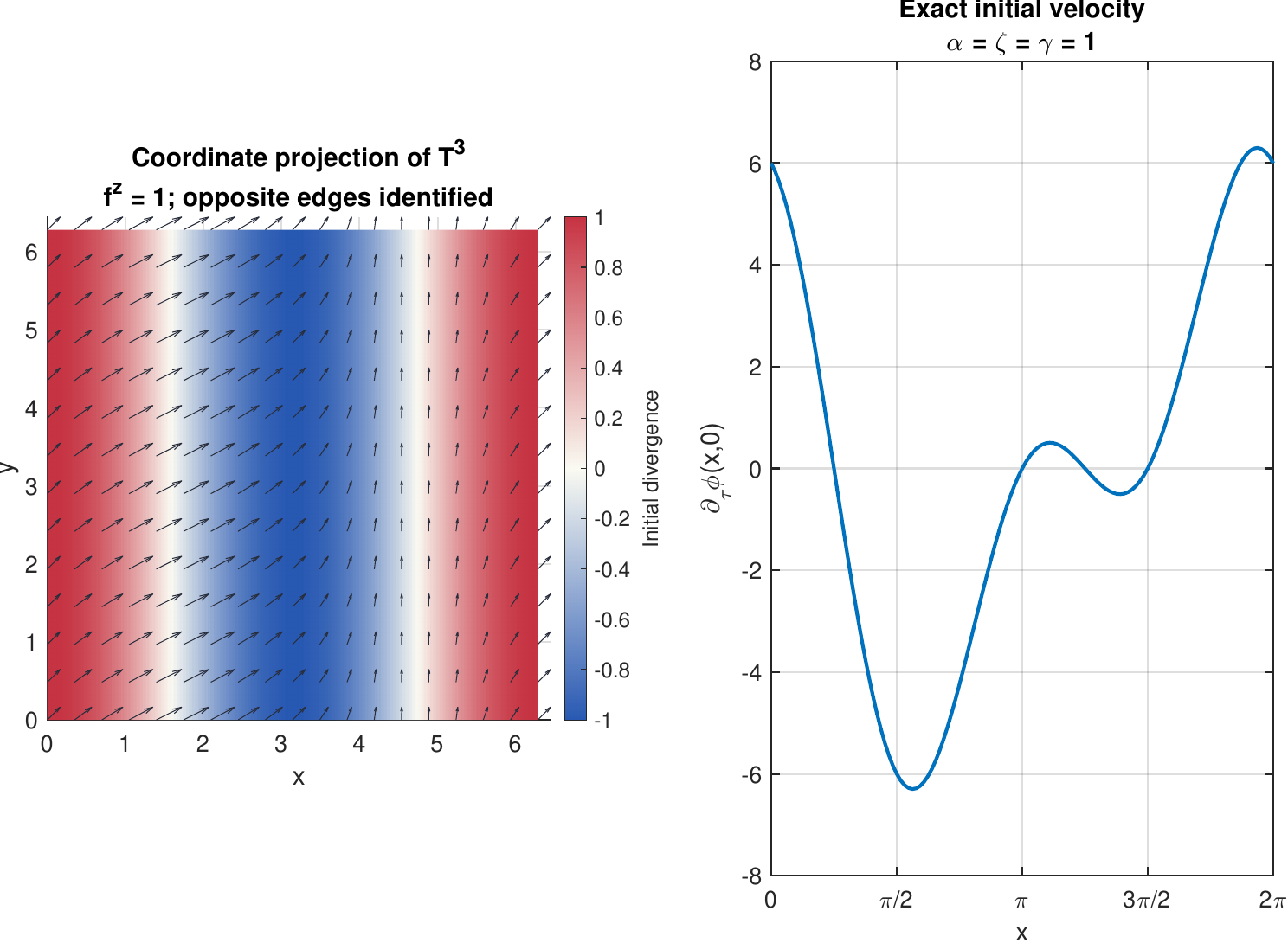}
\caption{Projection of $f$ onto $(x,y)$ and exact initial conformal velocity. Opposite sides are identified, the third field component is one. Color denotes the original divergence.}\label{Ex2_fig1}
\end{figure}
\begin{figure}[tbp]\centering
\manuscriptfigure{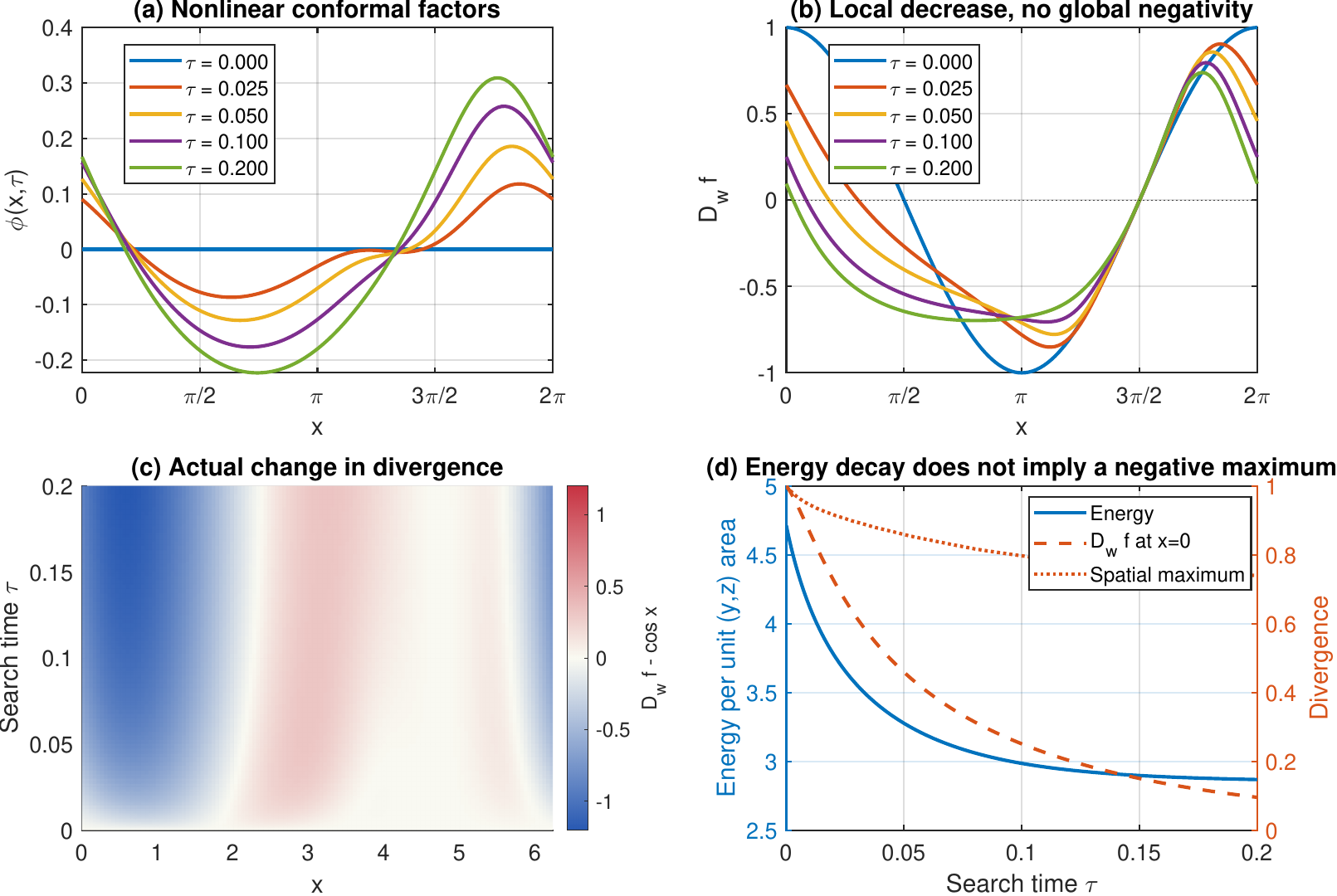}
\caption{Nonlinear unnormalized D-Yamabe evolution: conformal factor, divergence, its change, and energy. Divergence decreases at $x=0$ while its spatial maximum remains positive.}\label{Ex2_fig2}
\end{figure}
\end{example}

%%%%%%%%%%%%%%%%%%%%%%%%%%%%%%%%%%%%%%%%%%%%%%%%%%%

\section{D-Ricci Flow and the DeTurck Gauge}
\label{sec:ricci}

For arbitrary metric evolution use the Ebin metric $\langle h,k\rangle_{L^2(g)}=\int_M\langle h,k\rangle_gdV_g$ on symmetric $2$-tensors.
For a fixed smooth background metric $\bar g$, set
$V^k(g)=g^{ij}(\Gamma_{ij}^k(g)-\bar\Gamma_{ij}^k)$.
The classical Ricci--DeTurck operator $-2\Ric(g)+\mathcal L_{V(g)}g$ has scalar Laplace principal part and is equivalent to classical Ricci flow by pullback on common intervals of smooth existence~\cite{DeTurck1983}. The background is kept fixed when the initial metric varies.

\begin{theorem}
\label{th:RF_gradient_corrected}
Let $M$ be a compact oriented manifold without boundary, let $f\in\mathfrak{X}(M)$ be a fixed smooth vector field, let $\Psi:\R\to\R$ be smooth and convex, and let $\mu>0$. Consider the dissipative functional
$\mathcal{J}(g) = \int_M \Psi\bigl(\diver_g f\bigr) \, dV_g$.
The gradient of $\mathcal J$ relative to the Ebin $L^2$ metric is
\begin{equation}
\label{eq:grad_J}
\nabla_{L^2}\mathcal{J}(g) = \mathcal{T}_f(g),
\end{equation}
where the symmetric $2$-tensor field $\mathcal{T}_f(g)$ is defined by
\begin{equation}
\label{eq:T_f_def_corr}
\mathcal{T}_f(g) = \frac12\Bigl[ \Psi(\diver_g f) - \Psi'(\diver_g f)\,\diver_g f - \Psi''(\diver_g f)\, f(\diver_g f) \Bigr] g .
\end{equation}
Here $f(\diver_g f)$ is the derivative of $\diver_g f$ along $f$.

Define the \textit{D-Ricci flow} by adding the dissipative correction $-\mu\,\nabla_{L^2}\mathcal{J}(g)$ and the standard DeTurck gauge term to classical Ricci flow $\partial_t g = -2\Ric(g)$. The equation is
\begin{equation}
\label{eq:RF_flow_corr}
\partial_t g= -2 \Ric(g) - \mu\,\mathcal{T}_f(g) + \mathcal{L}_{V(g)}g .
\end{equation}
In particular, omitting the gauge term gives
\begin{equation}
\label{eq:RF_flow_reduced}
\partial_t g= -2\,\Ric(g) - \mu\,\mathcal{T}_f(g).
\end{equation}

For fixed $f$ and a smooth positive-definite initial metric,~\eqref{eq:RF_flow_corr} is strongly parabolic and has a unique local smooth solution. The same local existence and uniqueness hold in little H\"older spaces with initial metrics in $h^{2+\beta}$, $0<\beta<1$.

Suppose $(g,f)$ is a smooth solution of
\begin{equation}\label{eq:resh_flow}
\partial_tg=-2\Ric(g)-\mu\mathcal T_{f(t)}(g)+\mathcal L_{V(g)}g,
\qquad \partial_tf=\mathcal L_{V(g)}f.
\end{equation}
Let $\Phi_t$ be the diffeomorphisms generated by $-V(g(t))$ from the identity. Then $\widehat g=\Phi_t^*g$ and $\widehat f=\Phi_t^*f$ satisfy
$\partial_t\widehat g=-2\Ric(\widehat g)-\mu\mathcal T_{\widehat f}(\widehat g), \qquad \partial_t\widehat f=0.$
This identity for smooth solutions does not itself prove local solvability of the coupled system or assert its strict parabolicity as a system for $(g,f)$.

\end{theorem}

\begin{proof}
Put $q=\diver_gf$ and $h=\delta g$. The volume variation is $\delta dV_g=\tfrac12\tr_gh\,dV_g$, and, because $f$ is fixed, $\delta q=\tfrac12 f(\tr_gh)$. Integration by parts yields
$D\mathcal J(g)[h]=\tfrac12\int_M[\Psi(q)-\diver_g(\Psi'(q)f)]\tr_gh\,dV_g$.
Using $\diver_g(\Psi'(q)f)=\Psi'(q)q+\Psi''(q)f(q)$ and $\tr_gh=\langle h,g\rangle_g$ proves~\eqref{eq:grad_J}--\eqref{eq:T_f_def_corr}. Adding $-\mu\mathcal T_f$ to Ricci flow gives~\eqref{eq:RF_flow_reduced}, adding the DeTurck term gives~\eqref{eq:RF_flow_corr}.

In the positive-symbol convention, its principal symbol is $|\xi|_g^2h+\tfrac\mu4\Psi''(q)(f\cdot\xi)^2(\tr_gh)g$, whose pairing with $h$ is at least $|\xi|_g^2|h|_g^2$ for convex $\Psi$. This follows by taking the principal part of $\delta q=\tfrac12 f(\tr_gh)$ and using the Laplace principal part of the Ricci--DeTurck operator. Strong parabolicity and standard quasilinear theory give local solvability~\cite{Lieberman,PrussSimonett2016}. The penalty is natural only under simultaneous transformation of $(g,f)$, as follows from the pullback identity below, fixed-field gauged and ungauged problems are distinct. Neither Ricci curvature nor the gauge term is the negative gradient of $\mathcal J$. Ungauged stationarity is $-2\Ric(\hat g)-\mu\mathcal T_f(\hat g)=0$, gauged stationarity additionally contains $\mathcal L_{V(\hat g)}\hat g$.
Simultaneous transformation of the metric and vector field gives $\Phi^*\mathcal T_f(g)=\mathcal T_{\Phi^*f}(\Phi^*g)$. The derivative formula for pullback by a time-dependent diffeomorphism yields
$\frac d{dt}\Phi_t^*g =\Phi_t^*(\partial_tg-\mathcal L_{V(g)}g)$, $\frac d{dt}\Phi_t^*f =\Phi_t^*(\partial_tf-\mathcal L_{V(g)}f)$.
Substituting the system proves both identities.
\end{proof}

Writing $\mathcal T_{f,\Psi}$ to show the penalty dependence, \eqref{eq:T_f_def_corr} gives $\mathcal T_{f,\Psi+as+c}=\mathcal T_{f,\Psi}+cg/2$. In particular, an affine penalty produces the pure scaling term $cg/2$, and replacing $s^2/2$ by $(s+\gamma)^2/2$ adds only $-\mu\gamma^2g/4$ to the metric evolution. Thus the quadratic shift does not impose the threshold $\diver_g f\le-\gamma$. For a nonlinear penalty the term containing $f(\diver_g f)$ describes the derivative along the field.

\begin{proposition}
\label{th:extend_RF_deTurck}
Assume the hypotheses of Theorem~\ref{th:RF_gradient_corrected}, and let $g(t)$ be the maximal solution on $[0,T_{\max})$. If $T_{\max}<\infty$ and, for some $\beta\in(0,1)$, the metrics $g(t)$ remain uniformly equivalent to the fixed background metric $\bar g$, with
$\sup_{t<T_{\max}}\|g(t)\|_{C^{2,\beta}(\bar g)}<\infty,$
then the solution extends beyond $T_{\max}$.
\end{proposition}
\begin{proof}
The bounds ensure uniform parabolicity and a bounded $C^0$ right-hand side, hence a unique $C^0$ limit $g_T$ as $t\uparrow T_{\max}$. For $0<\beta'<\beta$, compactness of $C^{2,\beta}\hookrightarrow C^{2,\beta'}$ and uniqueness of that limit give convergence of the entire trajectory in $C^{2,\beta'}$. The limit remains positive definite. Restarting the local problem at $g_T$ with exponent $\beta'$ extends the solution.
\end{proof}

\begin{proposition}
\label{th:metric_spectral_convergence}
Let $M$ be compact and without boundary, assume the hypotheses of Theorem~\ref{th:RF_gradient_corrected}, and let $g_*$ be a smooth stationary state $F_R(g)=-2\Ric(g)-\mu\mathcal T_f(g)+\mathcal L_{V(g)}g$. Suppose that the linearization $L=DF_R(g_*)$, acting from the little H\"older space $h^{2+\beta}$ to $h^\beta$, $0<\beta<1$, has spectrum
\begin{equation}\label{eq:strict_spectral_gap}
\sup\{\operatorname{Re}z:z\in\sigma(L)\}\le-\omega<0.
\end{equation}
Then, for every $0<\omega'<\omega$, there exist $\delta,C>0$ such that, whenever $\|g_0-g_*\|_{h^{2+\beta}}<\delta$, the corresponding solution exists for all $t\ge0$ and
$\|g(t)-g_*\|_{h^{2+\beta}} \le C e^{-\omega't}\|g_0-g_*\|_{h^{2+\beta}}.$
If the stationary state satisfies $\operatorname{div}_{g_*}f
\le-2\gamma$ on a positively invariant domain $U$ for some $\gamma>0$, then sufficiently close initial data yield $\operatorname{div}_{g(t)}f\le-\gamma$ on $U$ for all $t\ge0$. Exponential contraction relative to the \emph{stationary} metric $g_*$ follows from the frozen-volume estimate~\eqref{eq:vol_exp_D}.
\end{proposition}
\begin{proof}
Theorem~\ref{th:RF_gradient_corrected} gives a uniformly parabolic quasilinear system near $g_*$. Its nonlinear operator acts smoothly from $h^{2+\beta}$ to $h^\beta$, and the linearization is sectorial with compact resolvent. Under~\eqref{eq:strict_spectral_gap}, the principle of linearized stability for quasilinear parabolic equations~\cite{PrussSimonett2016} gives the stated global solution and exponential estimate. Divergence depends continuously on the metric in $C^1$. The strict margin $2\gamma$ and smallness of $\delta$ therefore preserve its sign throughout $U$. The last conclusion concerns the fixed limiting metric. For volume relative to an evolving metric, the additional term $\frac12\operatorname{tr}_{g(t)}\dot g$ in~\eqref{eq:general_transport_th} must be included.
\end{proof}

\begin{proposition}
\label{cor:hyperbolic_stable_density}
Let $(M,\bar g)$ be a compact hyperbolic manifold of dimension $n\ge3$ and sectional curvature $-1$, with DeTurck background metric $\bar g$. Fix $\mu>0$ and
$\Psi(s)=\frac{4(n-1)}\mu+\frac{s^2}{2}$.
Let $x_*\in M$, and let $F$ be a smooth vector field that takes the form $F=-\sum_i x^i\partial_{x^i}$ on a coordinate ball $U$ near $x_*$. For sufficiently small $\delta>0$ and $f_\delta=\delta F$, the gauge-fixed $\mathcal D$-Ricci flow has a smooth stationary state $g_\delta$ close to $\bar g$ such that
$\operatorname{div}_{g_\delta}f_\delta<0 \quad\text{on a smaller ball }U_*\Subset U.$
All sufficiently close initial metrics yield global solutions converging exponentially to $g_\delta$. The ball $U_*$ may be chosen positively invariant under $f_\delta$, so the limiting metric provides local exponential phase-volume contraction.
\end{proposition}
\begin{proof}
At $\delta=0$, $\mathcal T_0(g)=2(n-1)g/\mu$ and $\bar g$ is stationary. For the resulting normalized flow
$\partial_tg=-2\operatorname{Ric}(g)-2(n-1)g +\mathcal L_{V(g)}g,$
with DeTurck background metric $\bar g$, the linearization at $\bar g$ acts on symmetric $2$-tensors by
\begin{equation}\label{eq:hyperbolic_L_explicit}
L_0h=\Delta_{\bar g}h+2h-2(\operatorname{tr}_{\bar g}h)\bar g,
\end{equation}
where $\bar\nabla$ is the connection of $\bar g$ and $\Delta_{\bar g}=\bar\nabla^k\bar\nabla_k$. The notation $S^2T^*M$ denotes the bundle of symmetric covariant second-order tensors. In particular, $\sup\sigma(L_0)\le-(n-2)<0$ on $L^2(S^2T^*M)$ and on the corresponding little H\"older spaces.

The Einstein linearization is $L_0=\Delta_L-2(n-1)$, which at sectional curvature $-1$ gives~\eqref{eq:hyperbolic_L_explicit}. The orthogonal splitting $h=h_0+u\bar g$, $\tr h_0=0$, is invariant: $L_0h_0=\Delta h_0+2h_0$ and $L_0(u\bar g)=(\Delta u-2(n-1)u)\bar g$. The Bochner--Koiso identity~\cite{Besse} reads $\|\bar\nabla h_0\|_2^2=\|\diver h_0\|_2^2+\tfrac12\|T(h_0)\|_2^2+n\|h_0\|_2^2$, where $T(h_0)_{ijk}=\bar\nabla_k(h_0)_{ij}-\bar\nabla_i(h_0)_{kj}$. Therefore $\langle L_0h_0,h_0\rangle_2\le-(n-2)\|h_0\|_2^2$, while the trace part has bound $-2(n-1)$. Self-adjointness and compact resolvent give the asserted gap. Elliptic regularity identifies the little H\"older spectrum with the same smooth eigentensors.

The strict gap makes the full linearization invertible, as required for the perturbation of the hyperbolic background~\cite{KroenckeEinstein2017}. The perturbation for $f_\delta=\delta F$ is $O(\delta^2)$, so the implicit function theorem gives $g_\delta=\bar g+O(\delta^2)$ and preserves the spectral gap. Proposition~\ref{th:metric_spectral_convergence} yields exponential stability. Since $\diver_{\bar g}F(x_*)=-n$, one has $\diver_{g_\delta}f_\delta=\delta\diver_{\bar g}F+O(\delta^3)<0$ on a smaller ball. There $F=-x$ points inward, ensuring invariance.
\end{proof}

\begin{proposition}
\label{th:global_DRicci_precompact}
Assume the hypotheses of Theorem~\ref{th:RF_gradient_corrected}, and let $g(t)$ be the maximal solution of~\eqref{eq:RF_flow_corr} on $[0,T_{\max})$. Suppose that, for some $\beta\in(0,1)$, there exist $C,c>0$ such that
\begin{equation}
\label{eq:DRicci_global_bounds}
 c^{-1}\bar g\le g(t)\le c\bar g,
 \qquad
 \sup_{t<T_{\max}}\|g(t)\|_{C^{2,\beta}(\bar g)}\le C.
\end{equation}
These bounds alone imply $T_{\max}=\infty$. Assume further that the global solution has finite action,
\begin{equation}
\label{eq:DRicci_finite_action}
 \int_0^\infty\|\partial_tg(t)\|_{L^2(\bar g)}^2dt<\infty.
\end{equation}
Then, for every $0<\beta_0<\beta$, the $\omega$-limit set of the trajectory in the little H\"older space $h^{2+\beta_0}$ is nonempty and consists of smooth stationary metrics of~\eqref{eq:RF_flow_corr}. If at least one metric $g_*\in\omega(g)$ is strictly spectrally stable, that is, for some $\omega>0$,
\begin{equation}
\label{eq:DRicci_omega_gap}
 \sup\{\operatorname{Re}z:z\in\sigma(D F_R(g_*))\}
 \le-\omega<0
\end{equation}
in the realization $h^{2+\beta_0}\to h^{\beta_0}$, then
$g(t) \rightarrow g_*$ in $C^k(M)$ for every $k\ge0$.
Moreover, for every $0<\omega'<\omega$, there exist $T$ and $C_{\omega'}>0$ such that
$\|g(t)-g_*\|_{h^{2+\beta_0}} \le C_{\omega'}e^{-\omega'(t-T)} \|g(T)-g_*\|_{h^{2+\beta_0}},\qquad t\ge T.$
\end{proposition}

\begin{proof}
The bounds~\eqref{eq:DRicci_global_bounds} and Proposition~\ref{th:extend_RF_deTurck} give global existence. Uniform parabolicity and translated-cylinder Schauder estimates give all higher bounds on positive-time tails, hence precompactness in $h^{2+\beta_0}$. The same estimates make $\partial_tg$ uniformly continuous in $L^2(\bar g)$ for $t\ge1$. Finite action then forces $\|\partial_tg(t)\|_2\to0$: otherwise disjoint intervals of uniformly positive action would contradict~\eqref{eq:DRicci_finite_action}.

For any convergent sequence $g(t_j)\to g_*$ in $h^{2+\beta_0}$, continuity of $F_R$ into $h^{\beta_0}$ and $F_R(g(t_j))=\partial_tg(t_j)\to0$ in $L^2$ give $F_R(g_*)=0$. Elliptic regularity makes $g_*$ smooth. If one limit point satisfies~\eqref{eq:DRicci_omega_gap}, the trajectory enters its stability neighborhood. Proposition~\ref{th:metric_spectral_convergence} controls the entire remaining trajectory with the stated exponential estimate. Smoothing yields convergence in every $C^k$.
The entry argument uses the local stability principle~\cite{PrussSimonettZacher2009}. When a spectrally stable stationary metric is already known to lie in $\omega(g)$, the same entry argument gives convergence without the finite-action assumption. Finite action is used above to prove stationarity of all limit points.
\end{proof}

If a bounded-below functional $\mathscr L$ along the orbit satisfies
$\partial_t\mathscr L(g(t))\le-c_0\|\partial_tg(t)\|_{L^2(\bar g)}^2$ for $c_0>0$, integration gives~\eqref{eq:DRicci_finite_action}, with the action bounded by $(\mathscr L(g(0))-\inf\mathscr L)/c_0$. The penalty term has its own gradient structure, and the Ricci and gauge terms affect its derivative. In particular, the monotonicity of Perelman's entropy~\cite{Perelman2002} does not transfer directly.
For $f=0$, $\Psi(0)=0$, and round $g_0$ of unit sectional curvature on $S^n$, $n\ge2$, the gauge-fixed flow with background $g_0$ has the solution $g(t)=(1-2(n-1)t)g_0$. The penalty vanishes, $V(g(t))=0$ because the scaling is spatially constant, and $\Ric(cg_0)=(n-1)g_0$ for $c>0$. The solution degenerates at $t=1/[2(n-1)]$.

\begin{example}
\label{Sec:Ricci_example}\label{e:riccitorus2}
Consider gauge-dependent departure from the conformal class.
On the same flat torus, take $f=\sin x\partial_x+\partial_y+\partial_z$ and $\Psi(s)=s^2/2$. Classical Ricci flow leaves $g_0$ fixed, with divergence $\cos x$. The torus $x=0$ consists of periodic orbits of period $2\pi$. For every metric $h$, write $dV_h=e^v dV_{g_0}$. Along each such orbit, $\diver_hf=1+(\partial_y+\partial_z)v$ has time average one. Thus no density, conformal metric, or general metric makes divergence negative on the entire invariant torus. More explicitly, averaging a density divergence over $(y,z)$ gives $\cos x+\sin x L'(x)\to1$ as $x\to0$, where $L$ is the mean logarithmic density. Negativity on a full strip $|x|<\epsilon$ is impossible.

Nevertheless, $\mathcal T_f(g_0)=\tfrac14(2-3\cos^2x)g_0$ and $g_t|_0=-\tfrac\mu4(2-3\cos^2x)g_0$, see Fig.~\ref{Ex3_fig1}. Write the initial conformal factor as $u=t u_1+O(t^2)$, where $u_1=-\mu(2-3\cos^2x)/8$ and $u_1'=-3\mu\sin x\cos x/4$. At $g=e^{2u(x)}g_0$, one has $\Ric_{xx}=-2u''$ and $\Ric_{yy}=\Ric_{zz}=-u''-(u')^2$. Consequently $\Ric_{xx}/g_{xx}-\Ric_{yy}/g_{yy}=(-u''+(u')^2)e^{-2u}=3\mu t\cos(2x)/4+O(t^2)$. The penalty tensor is conformal, so the ungauged right-hand side has a nonzero trace-free term of order $t$. Any smooth ungauged solution therefore leaves $[g_0]$ at order $t^2$.

For the gauge-fixed equation with background $g_0$, $V=-e^{-2u}\nabla u$ and $(\mathcal L_Vg)_{ij}=-2u_{ij}+4u_iu_j-2\delta_{ij}|\nabla u|^2$. Hence $-2\Ric(g)+\mathcal L_Vg=2u''g_0+2(u')^2dx^2$. The terms through second order in the metric are conformal, integrating the first trace-free right-hand side gives
\begin{equation}\label{eq:example_ricci_cubic}
 (g(\tau))^{\mathrm{TF}_{g_0}}
 =\frac23\tau^3(u_1')^2\left(dx^2-\frac13g_0\right)+O(\tau^4).
\end{equation}
Thus $g_{xx}-g_{yy}=\tfrac23\tau^3(u_1')^2+O(\tau^4)$. With the same prescribed $f$ held fixed in both equations, the first nonconformal term is cubic for the gauge-fixed evolution and quadratic for the ungauged evolution, conditional on smooth solvability of the latter. Formula~\eqref{eq:example_ricci_cubic} describes the initial-time expansion.

The exact symmetry-preserved diagonal class is $g=e^{2p(x,\tau)}dx^2+e^{2q(x,\tau)}(dy^2+dz^2)$, initially $p=q=0$. For $\mu=1$, set $D=\cos x+\sin x(p_x+2q_x)$ and $F=D^2/8+\sin x\,D_x/4$. Figures~\ref{Ex3_fig2} and~\ref{fig:ricci_obstruction} use the exact reduction~\eqref{eq:example_ricci_reduced},
\begin{equation}\label{eq:example_ricci_reduced}
 \begin{aligned}
 p_\tau=e^{-2p}(p_{xx}-p_x^2+2q_x^2)+F,~~~
 q_\tau=e^{-2p}q_{xx}+F.
 \end{aligned}
\end{equation}
Indeed, $\Ric_{xx}=-2(q_{xx}+q_x^2-p_xq_x)$, $\Ric_{yy}=\Ric_{zz}=-e^{2q-2p}(q_{xx}+2q_x^2-p_xq_x)$, and $V^x=e^{-2p}(p_x-2q_x)$. The volume density is $e^{p+2q}$, and $D(0,\tau)=1$ exactly. On $[0,0.2]$, the D-Yamabe grids and solver settings were reused, with $\left|\int D e^{p+2q}\,dx\right|$ monitored.

The ungauged equation has no stationary metric. Stationarity forces $\Ric(h)=c(x)h$, since $\mathcal T_f(h)$ is scalar. Schur's lemma makes $c$ constant. An Einstein metric on $\mathbb T^3$ must be flat: $c>0$ contradicts Bonnet--Myers, while $c<0$ gives constant negative sectional curvature, incompatible with $\pi_1(\mathbb T^3)=\mathbb Z^3$. Thus $\mathcal T_f(h)=0$, or $f(q_*)+q_*^2/2=0$ for $q_*=\diver_hf$. Integration along a periodic orbit on $x=0$ gives $q_*=0$, contradicting its mean one. This argument does not exclude stationary solutions of the gauge-fixed equation, which contains $\mathcal L_Vh$. Additional tensorial freedom therefore changes the evolution but does not remove the invariant divergence obstruction.
\begin{figure}[tbp]\centering
\manuscriptfigure{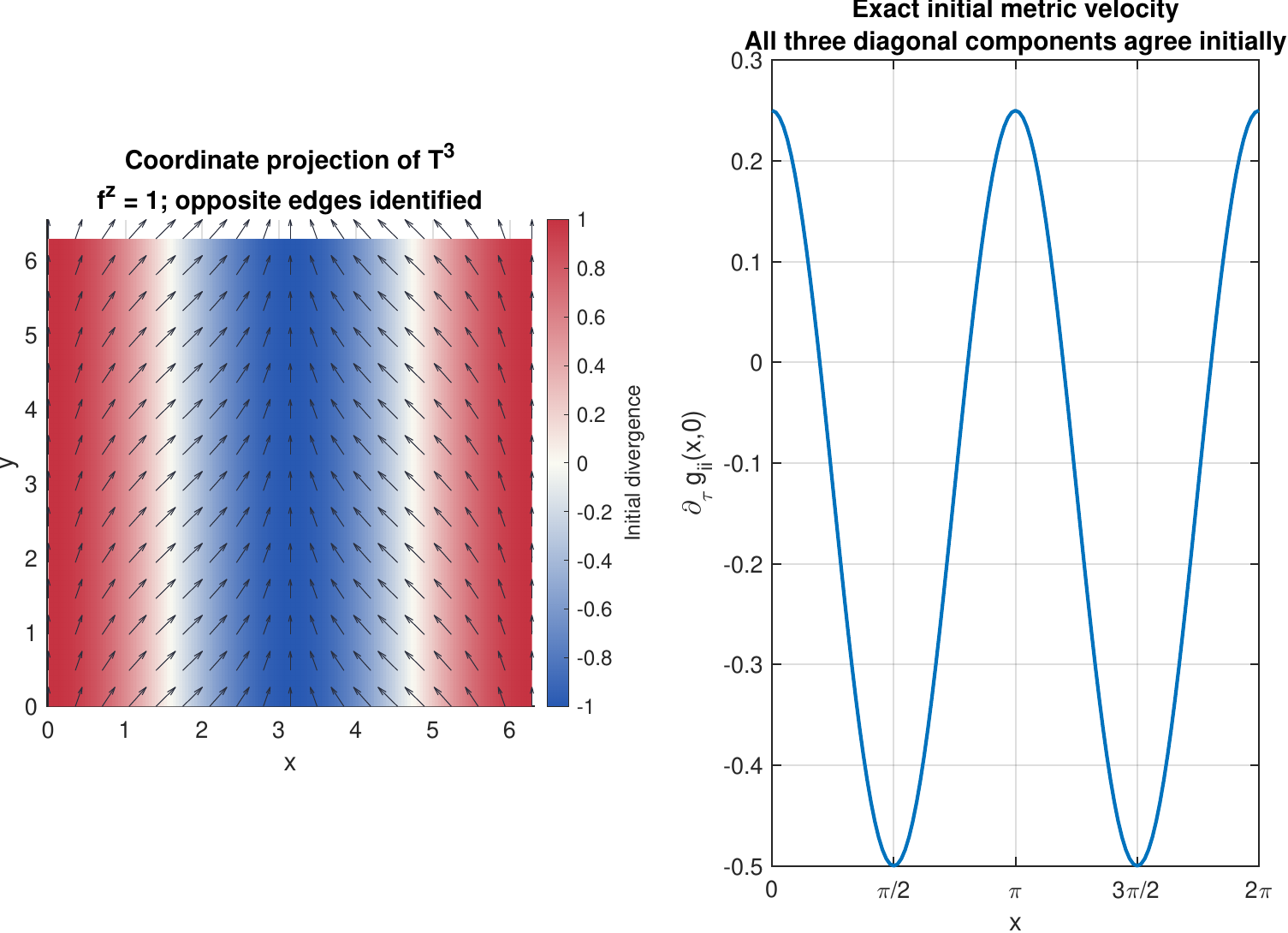}
\caption{Projection of the Ricci-example field and equal initial velocities of the diagonal metric components. Color denotes $\cos x$, the third field component is one.}\label{Ex3_fig1}
\end{figure}
\begin{figure}[tbp]\centering
\manuscriptfigure{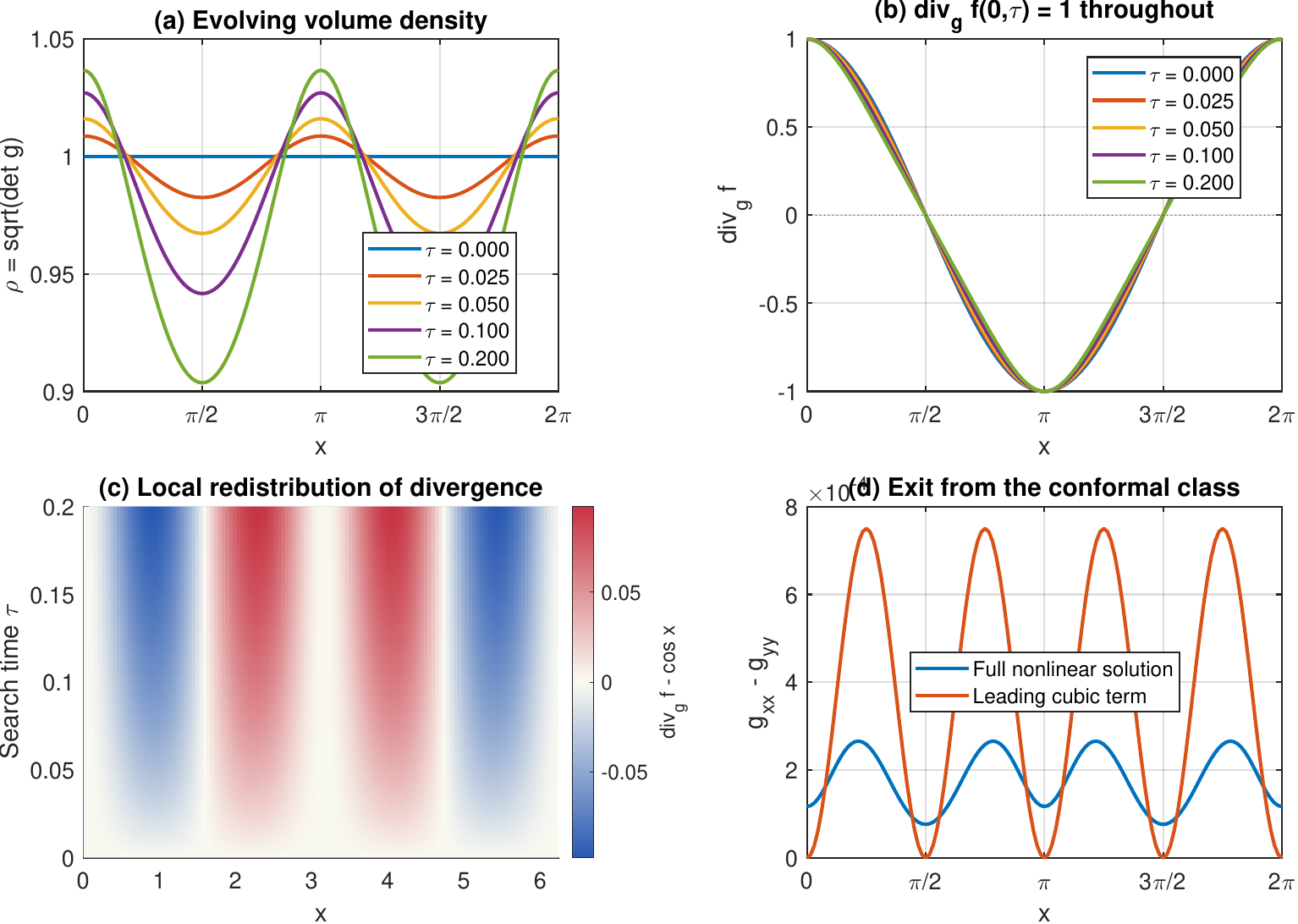}
\caption{Gauge-fixed D-Ricci evolution for fixed $f$: volume density, divergence, its change, and anisotropy compared with its cubic term. The two metric factors evolve independently, $D(0,\tau)=1$.}\label{Ex3_fig2}
\end{figure}
\begin{figure}[tbp]\centering
\manuscriptfigure{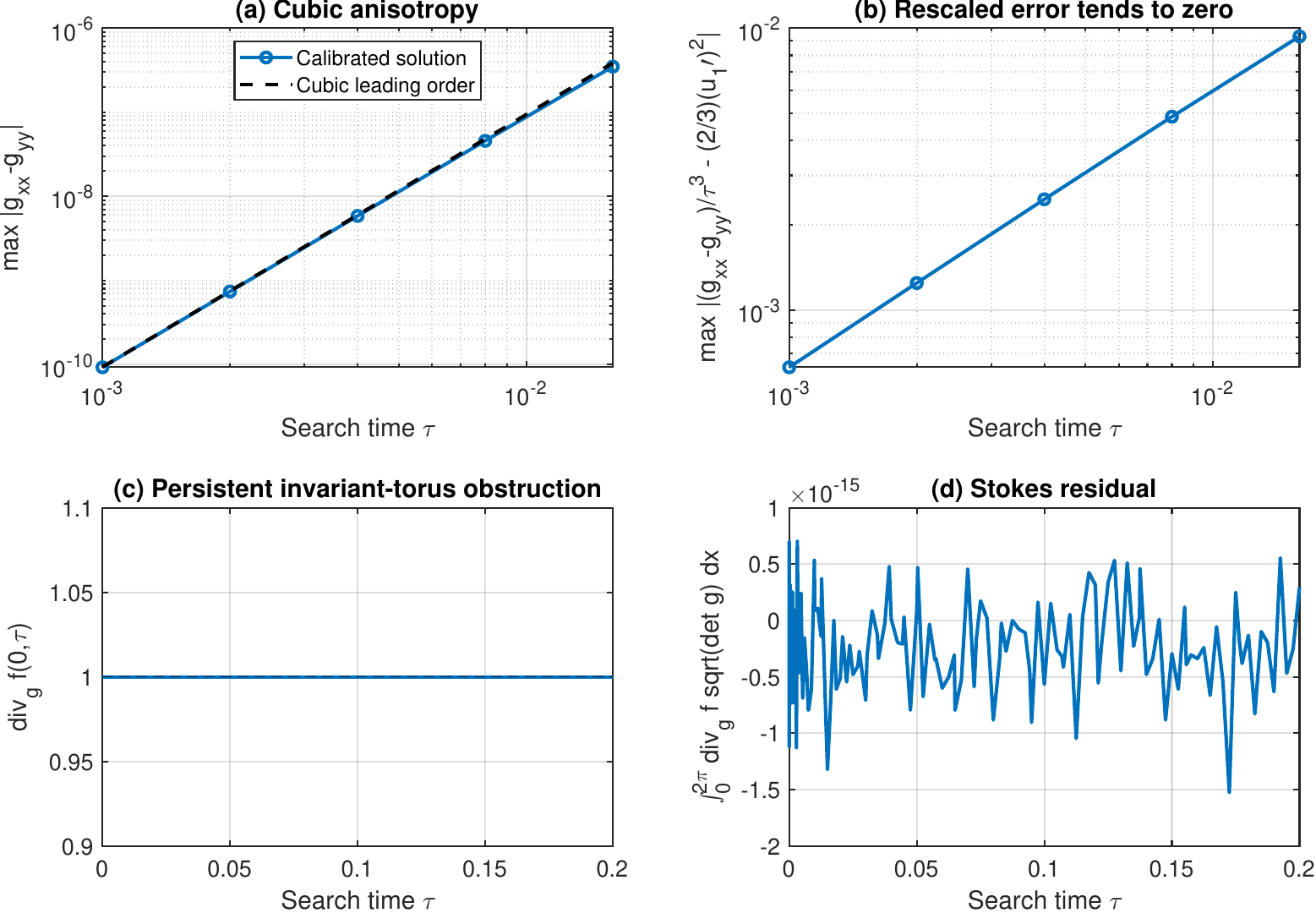}
\caption{Short-time cubic anisotropy and rescaled error, invariant-torus divergence, and Stokes residual. Departure from the conformal class preserves the invariant obstruction.}\label{fig:ricci_obstruction}
\end{figure}
\end{example}

%%%%%%%%%%%%%%%%%%%%%%%%%%%%%%%%%%%%%%%%%%%%%%%%%%%

\section{Stationary Limits and Physical Transport}
\label{sec:transport_limits}

\begin{proposition}
\label{th:limit_states}
Assume conditions ensuring local well-posedness of the corresponding flow.

1) If $\alpha\ne0$ and a solution $\psi(t)$ of the density flow~\eqref{eq:dens_psi} exists for every $t\ge 0$, converges to $\psi_\infty$ in the $C^2$ norm as $t\to\infty$, then $\psi_\infty$ satisfies the stationary equation
$\diver_g\bigl(\Phi'(D_{w_\infty^\alpha}f)\, f\bigr)=0$,
where $w_\infty=e^{\psi_\infty}$.

2) If a solution $\phi(t)$ of the D-Yamabe flow~\eqref{eq:YF_exact_corr} exists for every $t\ge 0$, converges to $\phi_\infty$ in the $C^2$ norm, then $\phi_\infty$ satisfies
$-(n-2)\alpha e^{\alpha n\phi_\infty}\widetilde R_{\alpha,\infty} + \zeta\alpha n\,\diver_g\!\bigl(\Phi'(D_{\hat w^\alpha}f)\,f\bigr)=0$,
where $\hat w=e^{n\phi_\infty}$.

3) If a solution $g(t)$ of the modified D-Ricci flow~\eqref{eq:RF_flow_corr} exists for every $t\ge 0$, converges in $C^2$ to a positive-definite metric $\hat g$, then $\hat g$ satisfies
$-2\,\Ric(\hat g) - \mu\,\mathcal{T}_f(\hat g) + \mathcal{L}_{V(\hat g)} \hat g = 0$,
where $V(\hat g)$ is the DeTurck vector field.

If the limiting divergence is uniformly negative on a positively invariant domain, then system~\eqref{eq:sys} is exponentially dissipative there relative to the constructed metric or weight.
\end{proposition}

\begin{proof}
Each right-hand side is continuous from $C^2$ states to $C^0$. In the density and conformal cases, $C^2$ convergence gives $C^1$ convergence of divergence, permitting passage through $\diver_g(\Phi'(D)f)$. Curvature and the DeTurck term converge in $C^0$ under $C^2$ metric convergence. For each search written as $z_t=F(z)$, this gives $z_t(t)\to F(z_\infty)$ in $C^0$. Since $z(t+1)-z(t)=\int_t^{t+1}F(z(s))\,ds\to0$, the limiting right-hand side is zero. Dividing by $\alpha$ in the density case requires $\alpha\ne0$. Freezing a limiting geometry with the stated sign gives~\eqref{eq:vol_exp_D}.
\end{proof}

\begin{corollary}\label{cor:frozen_search_limit}
Suppose a geometric search converges in $C^2$ to a smooth positive weight or a smooth Riemannian metric. On a positively invariant domain $U$, let its limiting volume form be
$\hat m=\hat w^\alpha dV_g$ for density evolution, $\hat m=e^{n\alpha\hat\phi}dV_g$ for conformal evolution, and $\hat m=dV_{\hat g}$ for general metric evolution.
If $\diver_{\hat m}f\le-\lambda_0<0$ on $U$, then, after freezing the geometry,
\begin{equation}\label{eq:frozen_search_volume}
 \int_{\varphi_t(D)}\hat m
 \le e^{-\lambda_0t}\int_D\hat m,
 \qquad D\subset U,\quad t\ge0.
\end{equation}
No convergence rate for the search enters this physical-time estimate.
\end{corollary}
\begin{proof}
Each limiting volume form is represented by a fixed Riemannian metric: respectively $\hat w^{2\alpha/n}g$, $e^{2\alpha\hat\phi}g$, and $\hat g$. The divergence identities in~\eqref{eq:volume_identity} and the frozen-volume estimate~\eqref{eq:vol_exp_D} give~\eqref{eq:frozen_search_volume}.
\end{proof}
An independent sign estimate is needed, for example the stopping test in Theorem~\ref{t:main2}. If the geometry is still changing during physical transport, the additional terms in~\eqref{eq:general_transport_th} remain necessary.

%%%%%%%%%%%%%%%%%%%%%%%%%%%%%%%%%%%%%%%%%%%%%%%%%%%

\section{Conclusions}
Controlled positive quadrature gives finite certification below the sharp rate, and vanishing regularization selects the minimum-norm feasible potential at a fixed margin. The annular examples verify volume decay with unit prefactor. The disk example shows that the sharp rate can lack a $C^1$ certificate.
The stable conformal and hyperbolic branches preserve the strict negativity of the background divergence on a small invariant ball. On the torus, a decrease in divergence or a nonconformal deformation can coexist with a contraction obstruction.

%%%%%%%%%%%%%%%%%%%%%%%%%%%%%%%%%%%%%%%%%%%%%%%%%%%

\appendix
\section{Convergence from a Gradient Inequality}
The normalized D-Yamabe proof uses the compactness and gradient-inequality argument \cite{Simon1983,Chill2003}. For a phase-field application of the same inequality \cite{WuGrasselliZheng2007}. Here the mean-zero projection removes the scaling direction.
\begin{lemma}
\label{th:abstract_LS}
Let $X$ be a Banach space densely and continuously embedded in a Hilbert space $H$, let $\mathcal E\in C^1(\mathcal O,\mathbb R)$, and let $\mathcal O\subset X$ be open. Suppose that $\mathcal G:\mathcal O\to H$ is continuous and is the $H$ gradient of $\mathcal E$ in the precise sense that
$D\mathcal E(u)[v]=\langle\mathcal G(u),v\rangle_H$ for $u\in\mathcal O$ and $v\in X$.
Consider a global strong solution
$u\in C([0,\infty),X)\cap C^1((0,\infty),H)$ of $\partial_tu=-\mathcal G(u)$,
whose energy $t\mapsto\mathcal E(u(t))$ is locally absolutely continuous and satisfies
$\frac d{dt}\mathcal E(u(t))=-\|\mathcal G(u(t))\|_H^2$ for almost every $t>0$.
Assume that the closure of $\{u(t):t\ge1\}$ is compact in $X$ and contained in $\mathcal O$. Let $\omega(u)$ be the set of limits in $X$ of sequences $u(t_j)$ with $t_j\to\infty$. Assume also that, near every critical point $v\in\omega(u)$, there exist $C_v>0$, $\theta_v\in(0,1/2]$, and an $X$ neighborhood $U_v$ such that
\begin{equation}\label{eq:LS_inequality_abstract}
 |\mathcal E(z)-\mathcal E(v)|^{1-\theta_v}
 \le C_v\|\mathcal G(z)\|_H,
 \qquad z\in U_v.
\end{equation}
Then there exists a critical point $u_\infty\in\omega(u)$ such that
$u(t)\to u_\infty$ in $X$ as $t\to\infty$.
If $\theta=1/2$ can be used in~\eqref{eq:LS_inequality_abstract} at $u_\infty$, convergence in $H$ is exponential. If $0<\theta<1/2$, then, for sufficiently large $t$,
$\mathcal E(u(t))-\mathcal E(u_\infty) \le C(1+t)^{-1/(1-2\theta)}$, $\|u(t)-u_\infty\|_H \le C(1+t)^{-\theta/(1-2\theta)}$.
\end{lemma}

\begin{proof}
Compactness and the energy identity give a finite energy limit $\mathcal E_\infty$ and square-integrable gradient. Along a sequence $t_j\to\infty$, compactness and continuity give $u(t_j)\to u_*$ in $X$, $\mathcal G(u_*)=0$, and $\mathcal E(u_*)=\mathcal E_\infty$. Set $e(t)=\mathcal E(u(t))-\mathcal E_\infty$. In an LS neighborhood of $u_*$, for $e>0$,
$-\partial_te^\theta=\theta e^{\theta-1}\|\mathcal G(u)\|_H^2\ge(\theta/C)\|\mathcal G(u)\|_H$.
If the trajectory first exited a fixed $X$ ball after times $t_j$, this inequality would make the corresponding $H$ displacement tend to zero. Precompactness would give an $X$-convergent subsequence of exit points with limit $u_*$, contradicting exit at the fixed radius. The trajectory therefore remains in one LS neighborhood, has finite $H$ length, and converges in $H$, precompactness upgrades this to $X$ convergence.
If $e$ reaches zero, the energy identity makes the subsequent solution stationary. Otherwise $e'\le-C^{-2}e^{2(1-\theta)}$ gives the stated energy decay, and $\|u(t)-u_\infty\|_H\le(C/\theta)e(t)^\theta$ gives the distance rate, exponential for $\theta=1/2$.
\end{proof}

%%%%%%%%%%%%%%%%%%%%%%%%%%%%%%%%%%%%%%%%%%%%%%%%%%%

\setlength{\bibsep}{0pt}


\begin{thebibliography}{99}

\bibitem{Rantzer01}
A. Rantzer, A dual to Lyapunov's stability theorem, Systems \& Control Letters 42 (3) (2001) 161--168. \href{https://doi.org/10.1016/S0167-6911(00)00087-6}{DOI}.

\bibitem{PrajnaParriloRantzer2004}
S. Prajna, P.A. Parrilo, A. Rantzer, Nonlinear control synthesis by convex optimization, IEEE Transactions on Automatic Control 49 (2) (2004) 310--314. \href{https://doi.org/10.1109/TAC.2003.823000}{DOI}.

\bibitem{Lohmiller98}
W. Lohmiller, J.-J.E. Slotine, On contraction analysis for non-linear systems, Automatica 34 (6) (1998) 683--696.

\bibitem{Forni14}
F. Forni, R. Sepulchre, A differential Lyapunov framework for contraction analysis, IEEE Transactions on Automatic Control 59 (3) (2014) 614--628.

\bibitem{Furtat26_Flow}
I.B. Furtat, Weighted phase volume stability: dissipativity and geometric interpretation, arXiv:2610.06045 (2026).

\bibitem{Furtat26_Density}
I.B. Furtat, Convex density design for regional phase-volume contraction of nonlinear systems, arXiv:2610.06066 (2026).

\bibitem{TobascoGoluskinDoering2018}
I. Tobasco, D. Goluskin, C.R. Doering, Optimal bounds and extremal trajectories for time averages in nonlinear dynamical systems, Physics Letters A 382 (6) (2018) 382--386. \href{https://doi.org/10.1016/j.physleta.2017.12.023}{DOI}.

\bibitem{Hamilton1982}
R.S. Hamilton, Three-manifolds with positive Ricci curvature, Journal of Differential Geometry 17 (2) (1982) 255--306. \href{https://doi.org/10.4310/jdg/1214436922}{DOI}.

\bibitem{Ye1994}
R. Ye, Global existence and convergence of Yamabe flow, Journal of Differential Geometry 39 (1) (1994) 35--50. \href{https://doi.org/10.4310/jdg/1214454674}{DOI}.

\bibitem{Brendle07}
S. Brendle, Convergence of the Yamabe flow for arbitrary initial energy, Journal of Differential Geometry 69 (2) (2005) 217--278. \href{https://doi.org/10.4310/jdg/1121449107}{DOI}.

\bibitem{Yan22}
Z. Yan, Convergence of the weighted Yamabe flow, Differential Geometry and its Applications 84 (2022) 101922. \href{https://doi.org/10.1016/j.difgeo.2022.101922}{DOI}.

\bibitem{Chen2026Willmore}
Y. Chen, Volume preserving Willmore flow in a generalized Cahn--Hilliard flow, J. Math. Anal. Appl. 558 (2) (2026) 130451. \href{https://doi.org/10.1016/j.jmaa.2026.130451}{DOI}.

\bibitem{GieslHafsteinMehrabinezhad2021}
P. Giesl, S. Hafstein, I. Mehrabinezhad, Computation and verification of contraction metrics for periodic orbits, J. Math. Anal. Appl. 503 (2) (2021) 125309. \href{https://doi.org/10.1016/j.jmaa.2021.125309}{DOI}.

\bibitem{WuGrasselliZheng2007}
H. Wu, M. Grasselli, S. Zheng, Convergence to equilibrium for a parabolic--hyperbolic phase-field system with dynamical boundary condition, J. Math. Anal. Appl. 329 (2) (2007) 948--976. \href{https://doi.org/10.1016/j.jmaa.2006.07.011}{DOI}.

\bibitem{FeehanMaridakis2017}
P.M.N. Feehan, M. Maridakis, \L{}ojasiewicz--Simon gradient inequalities for analytic and Morse--Bott functions on Banach spaces, Journal f\"ur die reine und angewandte Mathematik 765 (2020) 35--67. \href{https://doi.org/10.1515/crelle-2019-0029}{DOI}.

\bibitem{DeTurck1983}
D.M. DeTurck, Deforming metrics in the direction of their Ricci tensors, Journal of Differential Geometry 18 (1) (1983) 157--162. \href{https://doi.org/10.4310/jdg/1214509286}{DOI}.

\bibitem{Liouville}
J. Liouville, Note sur la théorie de la variation des constantes arbitraires, J. Math. Pures Appl 3 (1838) 342--349.

\bibitem{Reynolds}
O. Reynolds, Papers on Mechanical and Physical Subjects. Cambridge: Cambridge University Press, 1903.

\bibitem{Furtat26_Arc}
I.B. Furtat, Weighted phase volume method in stability analysis: integral criteria and ellipsoidal reachable sets, arXiv:2607.05033 (2026).

\bibitem{Kunita1990}
H. Kunita, Stochastic Flows and Stochastic Differential Equations. Cambridge Studies in Advanced Mathematics. Vol. 24. Cambridge: Cambridge University Press, 1990.

\bibitem{Brezis}
H. Brezis, Opérateurs maximaux monotones et semi-groupes de contractions dans les espaces de Hilbert, North-Holland Mathematics Studies, vol.~5, North-Holland, Amsterdam, 1973.

\bibitem{Besse}
A.L. Besse, Einstein manifolds. Berlin: Springer, 1987. \href{https://doi.org/10.1007/978-3-540-74311-6}{DOI}.

\bibitem{PrussSimonett2016}
J. Pr\"uss, G. Simonett, Moving Interfaces and Quasilinear Parabolic Evolution Equations. Monographs in Mathematics, vol.~105. Birkh\"auser, 2016. \href{https://doi.org/10.1007/978-3-319-27698-4}{DOI}.

%\bibitem{Ladyzhenskaya}
%O.A. Ladyzhenskaya, V.A. Solonnikov, N.N. Ural'tseva, Linear and Quasilinear Equations of Parabolic Type (in Russian). Moscow: Nauka, 1967.

\bibitem{Lieberman}
G.M. Lieberman, Second Order Parabolic Differential Equations. World Scientific, 1996.

\bibitem{KroenckeEinstein2017}
K. Kr\"oncke, Stability of Einstein metrics under Ricci flow, Communications in Analysis and Geometry 28 (2) (2020) 351--394. \href{https://doi.org/10.4310/CAG.2020.v28.n2.a5}{DOI}.

\bibitem{PrussSimonettZacher2009}
J. Pr\"uss, G. Simonett, R. Zacher, On convergence of solutions to equilibria for quasilinear parabolic problems, Journal of Differential Equations 246 (10) (2009) 3902--3931. \href{https://doi.org/10.1016/j.jde.2008.10.034}{DOI}.

\bibitem{Perelman2002}
G. Perelman, The entropy formula for the Ricci flow and its geometric applications, arXiv:math/0211159 (2002).

\bibitem{Simon1983}
L. Simon, Asymptotics for a class of non-linear evolution equations, with applications to geometric problems, Annals of Mathematics 118 (3) (1983) 525--571. \href{https://doi.org/10.2307/2006981}{DOI}.

\bibitem{Chill2003}
R. Chill, On the \L{}ojasiewicz--Simon gradient inequality, Journal of Functional Analysis 201 (2) (2003) 572--601. \href{https://doi.org/10.1016/S0022-1236(02)00102-7}{DOI}.

\end{thebibliography}
\end{document}